\documentclass[aos,MSNbibl,seceqn,dvips]{arximspdf}
\usepackage{graphicx}
\doi{10.1214/14-AOS1248} 
\volume{42}
\issue{5}
\pubyear{2014}
\firstpage{1970}
\lastpage{2002}
\docsubty{FLA}

\makeatletter
\newcommand{\Exp}{\operatorname{exp}}
\newcommand{\rrvert}{\vert}
\newcommand{\rrVert}{\Vert}
\newcommand{\llvert}{\vert}
\newcommand{\llVert}{\Vert}
\newtheorem{prop}{Proposition}[section]
\newtheorem{teo}{Theorem}[section]
\newtheorem{lem}{Lemma}[section]
\newproclaim{rem}{Remark}
\newproclaim{ass}{Assumption}[section]
\newproclaim{defn}{Definition}
\newcommand{\ind}{\mathbf{1}}
\newcommand{\N}{\mathbb{N}}
\newcommand{\OO}{\mathcal{O}}
\newcommand{\NN}{\mathcal{N}}
\newcommand{\aexp}{\mathfrak{a}}
\newcommand{\z}{\mathfrak{z}}
\newcommand{\cf}{\mathfrak{c}}
\newcommand{\bd}{\mathfrak{b}}
\newcommand{\zz}{x}
\newcommand{\mm}{m}
\makeatother

\begin{document}
\begin{frontmatter}

\title{Adaptive function estimation in nonparametric regression with
one-sided errors}
\runtitle{Adaptation with one-sided errors}

\begin{aug}
\author[1]{\fnms{Moritz}~\snm{Jirak}\thanksref{t1}\ead[label=e1]{jirak@math.hu-berlin.de}},
\author[2]{\fnms{Alexander}~\snm{Meister}\corref{}\ead[label=e2]{alexander.meister@uni-rostock.de}}
\and
\author[1]{\fnms{Markus}~\snm{Rei{\ss}}\ead[label=e3]{mreiss@mathematik.hu-berlin.de}}
\runauthor{M. Jirak, A. Meister and M. Rei{\ss}}
\affiliation{Humboldt Universit\"{a}t zu Berlin,
Universit\"at Rostock
and\\
Humboldt Universit\"{a}t zu Berlin}
\address[1]{M. Jirak\\
M. Rei{\ss}\\
Institut f\"{u}r Mathematik\\
Humboldt-Universit\"at zu Berlin\\
Unter den Linden 6\\
D-10099 Berlin\\
Germany\\
\printead{e1}\\
\phantom{E-mail: }\printead*{e3}}
\address[2]{A. Meister\\
Institut f\"{u}r Mathematik\\
Universit\"at Rostock\\
D-18051 Rostock\\
Germany\\
\printead{e2}}
\end{aug}
\thankstext{t1}{Supported by Deutsche Forschungsgemeinschaft via FOR1735
{Structural Inference in Statistics: Adaptation and Efficiency}.} 

\received{\smonth{11} \syear{2013}}
\revised{\smonth{4} \syear{2014}}

%
\begin{abstract}
We consider the model of nonregular nonparametric regression where
smoothness constraints are imposed on the regression function $f$ and
the regression errors are assumed to decay with some sharpness level at
their endpoints. The aim of this paper is to construct an adaptive
estimator for the \mbox{regression} function $f$. In contrast to the standard
model where local averaging is fruitful, the nonregular conditions
require a substantial different treatment based on local extreme
values. We study this model under the realistic setting in which both
the smoothness degree $\beta> 0$ and the sharpness degree $\mathfrak
{a}\in(0, \infty)$ are unknown in advance. We construct adaptation
procedures applying a nested version of Lepski's method and the
negative Hill estimator which show no loss in the convergence rates
with respect to the general \mbox{$L_q$-}risk and a logarithmic loss with
respect to the pointwise risk. Optimality of these rates is proved for
$\mathfrak{a}\in(0, \infty)$. Some numerical simulations and an
application to real data are provided.
\end{abstract}

%
\end{frontmatter}

\section{Introduction} \label{1}

In the standard model of nonparametric regression, the data
%
\begin{equation}
\label{eq11} Y_j = f(x_j) + \varepsilon_j,
\qquad j=1,\ldots,n
\end{equation}
are observed. In this paper, in contrast to classical theory, the
observation errors $(\varepsilon_j)$ are not assumed to be centred,
but to have certain support properties. This is motivated from many
applications where rather the support than the mean properties of the
noise are known and where the regression function $f$ describes some
frontier or boundary curve. Below we shall discuss concrete
applications to sunspot data and annual sport records. Typical
economical examples include auctions where the bidders' private values
are inferred from observed bids (see Guerre et al.~\cite{guerre2000}
or Donald and Paarsch~\cite{paarsch1993}) and note the extension to
bid and ask prices in financial data. Related phenomena arise in the
context of inference for deterministic production frontiers, where it
is assumed that $f$ is concave (convex) or monotone.

A pioneering contribution in this area is due to Farrell~\cite
{farrell1957}, who introduced data envelopment analysis (DEA), based on
either the conical hull or the convex hull of the data. This was
further extended by Deprins et al.~\cite{deprins1984} to the free
disposal Hull (FDH) estimator, whose properties have been extensively
discussed in the literature; see, for instance, Banker~\cite
{banker1993}, Korostelev et al.~\cite{korostelev1995a}, Kneip et
al.~\cite{Kneip1998,kneip2008}, Gijbels et al.~\cite{gijbels1999}, Park
et al.~\cite{park2000,park2010}, Jeong and Park~\cite{jeong2006}
and Daouia et al.~\cite{daouiasimar2010}. The issue of stochastic
frontier estimation goes back to the works of Aigner et al.~\cite
{Aigner1977} and Meeusen and van den Broeck~\cite{meeusen1977}; see
also the more recent contributions of Kumbhakar et al.~\cite
{Kumbhakar2007}, Park et al.~\cite{Park2007} and Kneip et al.~\cite
{kneip2012}.

In a general nonparametric setting the accuracy of the
estimator heavily depends on the average number of observations in the
vicinity of the support boundary. The key quantity is the
sharpness $\aexp_{\zz}>0$ of the distribution function $F_{\zz}$ of
$\varepsilon_j$ at $\zz=x_j$, which in its simplest case has
polynomial tails
%
\begin{equation}
\label{defstandard} F_{\zz}(y) = 1 - \cf_{\zz}'
\llvert y\rrvert^{\aexp_{\zz}} + \OO\bigl(\llvert y\rrvert^{\aexp
_{\zz} + \delta}
\bigr)\qquad\mbox{with }\cf_{\zz}', \delta> 0\mbox{ as }y
\to0.
\end{equation}
The cases $0 < \aexp_{\zz} < 1$, $\aexp_{\zz} = 1$ and $\aexp_{\zz
} > 1$ are sometimes called \emph{sharp boundary}, \emph{fault-type
boundary} and \emph{nonsharp boundary}. From a theoretical \mbox{perspective}
noise models with $\aexp_{\zz}\in(0,2)$ are nonregular (e.g.,
Ibragimov and Hasminskii~\cite{ibrakhasminskii1981}) since they
exhibit nonstandard statistical theory already in the parametric case.
Chernozhukov and Hong~\cite{chernozhukovhong2004} discuss extensively
parametric efficiency of maximum-likelihood and Bayes estimators in
this context and show their relevance in \mbox{economics}.

From a nonparametric statistics point of view, Korostelev and
Tsybakov~\cite{korostelevtsyb1993} and Goldenshluger and Zeevi~\cite
{goldenshluger2006}
treat a variety of boundary estimation problems. The focus is on
applications in image
recovery and is mathematically\vspace*{1pt} and practically substantially different
from ours. The optimal convergence rate $n^{(-2 \beta)/(\aexp\beta+1)}$ over $\beta$-H\"older classes of regression functions $f$
depends heavily on $\aexp$ (not assumed to be varying in $x$); for
$\aexp_{\zz}\in(0,2)$ it is faster than for local averaging
estimators in standard mean regression and can even become faster than
the regular squared parametric rate $n^{-1}$. Hall and van
Keilegom~\cite{hallkeilegom2009} study a local-linear estimator in a closer
related nonparametric regression model and establish minimax optimal
rates in $L_2$-loss if the smoothness and sharpness parameters $\beta
\in(0,2]$ and $\aexp>0$ are known. Earlier contributions in a related
setup are due to H\"{a}rdle et al.~\cite{haerdle1995}, Hall et al.~\cite
{hall1997,hall1998}
and Gijbels and Peng~\cite{gijbels2000}. If the support of
$(\varepsilon_j)$ is not one-sided, but symmetric like $[-a,a]$ and
$\beta\le1$, $\aexp=1$, M\"uller and Wefelmeyer~\cite
{muellerwefelmayer2010} have shown that mid-range estimators attain
also these better rates.
Recently, Meister and Reiss~\cite{meisterreiss2013} have proved
strong asymptotic equivalence in Le Cam's sense between a nonregular
nonparametric regression model for $\aexp=1$ and a continuous-time
Poisson point process experiment.

All the references above consider a theoretically optimal bandwidth
choice which depends on the unknown quantities $\aexp$ and/or
$\beta$. Completely data-driven adaptive procedures have been rarely
considered in the literature because the \mbox{intrinsically} nonlinear
inference and the nonmonotonicity of the stochastic and approximation
error terms block popular concepts from mean regression like
cross-validation or general unbiased risk estimation; cf. the
discussion in Hall and Park~\cite{hall2004}.
Recently, Chichignoud~\cite{chichi2012} was able to produce a $\beta
$-adaptive minimax optimal estimator, which, however, uses a Bayesian
approach hinging on the assumption that the law of the errors
$(\varepsilon_j)$ is perfectly known in advance (in fact, after log
transform a uniform law is assumed). Moreover, a log factor due to
adaptation is paid, which is natural only under pointwise loss. It
remained open whether under a global loss function like an $L_q$-norm
loss adaptation without paying a log factor is possible. For regular
nonparametric problems Goldenshluger and Lepski~\cite
{goldenshlugerlepski2011} study adaptive methods and convergence rates
with respect to general $L_q$-loss which is much more involved in the
general case $q\ge1$ than for $q=2$.

It is therefore of high interest, both from a theoretical and a
practical perspective, to establish a fully data-driven estimation
procedure where the error distribution and the regularity of the
regression function are unknown and to analyze it under local
(pointwise) and global ($L_q$-norm) loss. In particular, neither $\aexp
$ nor $\beta$ that determine the optimal convergence rate are fixed in
advance. In this paper we introduce a fully data-driven ($\aexp,\beta
$)-adaptive procedure for estimating $f$ and prove that it is
minimax optimal over $\aexp, \beta> 0$.

To ease the presentation, we restrict to equidistant design points $x_j
= j/n$ on $[0,1]$ and regression errors $(\varepsilon_j)$ which are
concentrated on the interval $(-\infty,0]$. Given $\zz\in[0,1]$ and
an open neighborhood $\NN(\zz) \subseteq[0,1]$, the function $f\dvtx \NN
(\zz) \to\mathbb{R}$ is supposed to lie in the H\"older class
$H_{\NN(\zz)}(\beta,L)$ with $\beta,L >0 $. Note that $\beta=
\beta_{\zz}$ and $L = L_{\zz}$ may vary in $\zz$. The $\langle
\beta\rangle$-derivatives of all $f \in H_{\NN(\zz)}(\beta,L)$ satisfy
\[
\bigl\llvert f^{(\langle\beta\rangle)}(y) - f^{(\langle\beta\rangle
)}(z)\bigr\rrvert\leq L
\llvert y-z\rrvert^{\beta- \langle\beta\rangle
},\qquad y,z \in\NN(\zz).
\]
Here $\langle\beta\rangle=\max\{ m\in\N_0\dvtx m<\beta\}$ is the
largest integer strictly smaller than $\beta$.

We consider the case where the $\varepsilon_j$ are independent with
individual distribution function $F_{x_j}$ and tail quantile function
\[
\mathcal{U}_{x_j}(y) = F_{x_j}^{\leftarrow} (1 - 1/y ),
\]
where\vspace*{1pt} $F_{x_j}^{\leftarrow}$ denotes the generalized inverse of
$F_{x_j}$. Weakening the polynomial tail behavior in (\ref
{defstandard}), our key structural condition is that for each $\zz\in
[0,1]$, there exist $\aexp_{\zz},\cf_{\zz} > 0$, $\bd_{\zz} \in
\mathbb{R}$ and a slowly varying function $l_{\zz}(y)$, such that
%
\begin{equation}
\label{defU} \mathcal{U}_{x}(y) = -\cf_x
y^{-1/\aexp_x} l_x(y),
\end{equation}
where $l_x(y)$ satisfies uniformly for $x \in[0,1]$ condition
%
\begin{equation}
\label{defL} l_x(y) = \log(y)^{\bd_x} + \mbox{\scriptsize$
\mathcal{O}$} \bigl(\log(y)^{\bd_x - 1} \bigr)\qquad\mbox{as }y \to
\infty.
\end{equation}
If (\ref{defstandard}) holds, then (\ref{defU}), (\ref{defL}) are
valid with $\bd_{\zz} = 0$
(note that $\cf_{\zz} \neq\cf_{\zz}'$ in general; see Lemma~\ref{lemquantcoomp} for the precise relation).
The polynomial tail condition (\ref{defstandard}) is one of the
standard models in the literature; see
de Haan and Ferreira~\cite{dehaanbook2006}, H\"{a}rdle et al.~\cite
{haerdle1995}, Hall and van Keilegom~\cite{hallkeilegom2009} or
Girard et al.~\cite{girard2013}. In this context, so called \textit{second order conditions} are inevitable whenever one is interested in
convergence rates or limit distributions involving estimates of $\aexp
_{\zz}$; see Beirlant et al.~\cite{beirlantteugelsbook2004}, de Haan
and Ferreira~\cite{dehaanbook2006} or Falk et al.~\cite{falkhr94}.
Our second order condition (\ref{defL}) is rather mild when compared
to examples from the literature; cf.~\cite
{beirlantteugelsbook2004,dehaanbook2006,falkhr94,girard2013,hallkeilegom2009,haerdle1995}.
As will be explained in Section~\ref{32}, a more general formulation
seems to be impossible.


Let us point out two main conceptual results of this paper. First, we
wish to extend the existing theory beyond the limitation $\beta_{\zz
}\le2$ imposed by locally constant or linear approximations and to
have a clear notion of stochastic and deterministic error for the
nonlinear estimators. To this end we develop a linear program in terms
of general local polynomials, based on a quasi-likelihood method,
because the definition in Hall and van Keilegom~\cite
{hallkeilegom2009} does not extend to polynomials of degree~2 or more
in our setup. Then Theorem~\ref{teocon} below yields for the estimator
a nontrivial decomposition in approximation and stochastic error. This
decomposition is a key result for our analysis, and permits us to
address the adaptation problem in full generality, thus abolishing the
blockade mentioned in Hall and Park~\cite{hall2004}.
We can consider not only pointwise, but also the global $L_q$-norm as
risk measure for the whole range $q \in[1, \infty)$. Technically, the
optimal $L_q$-adaptation is much more demanding compared to the
pointwise risk. It requires very tight deviation bounds since no
additional $\log n$-factor widens the margin.

For adaptive bandwidth selection, we apply a nested variant of
the\break
Lepski~\cite{lepski1990} procedure with pre-estimated critical
values. Careful adaptive pre-estimation is necessary since the
distribution of $(\varepsilon_j)$ is unknown and allowed to vary in
$\zz$. The fact that the underlying sample $(Y_j)$ is inhomogeneously
shifted by $f$ adds another level of complexity for the estimation\break of
$\aexp_{\zz}$ and $\bd_{\zz}$, which needs to be addressed by
translation invariant estimators. The remarkable result of Theorem~\ref{TL2upperbound} is
that for general $L_q$-loss\break we obtain the rate
$n^{(-2 \beta)/(\aexp\beta+1)} (\log n)^{(2\aexp_{\zz}
\bd_{\zz} \beta_{\zz})/(\aexp_{\zz} \beta_{\zz}+1)}$ of
convergence, the same as in the case of known (global) H\"older
regularity $\beta$ and known\break distribution of $(\varepsilon_j)$. For
pointwise loss the rate deteriorates to\break $(n /\log n)^{(-2 \beta
_{\zz})/({\aexp_{\zz} \beta_{\zz}+1})} (\log n)^{(2\aexp_{\zz
} \bd_{\zz} \beta_{\zz})/(\aexp_{\zz} \beta_{\zz}+1)}$; see
Theorem~\ref{Tpointwise} below.

In Section~\ref{seclowerbounds} it is shown that all our rates are
minimax optimal for adaptive estimation. For regular mean regression
these rates, inserting $\aexp_{\zz}=2$ and $\bd_{\zz} = 0$, and
particularly the payment for adaptation on $\beta_{\zz}$ under
pointwise loss are well known. A priori it is, however, not at all
obvious that in the nonregular case with Poisson limit experiments
(Meister and Reiss~\cite{meisterreiss2013}) exactly the same
factor appears. Interestingly, we do not pay in the convergence rates
for not knowing $\aexp_{\zz}$, $\bd_{\zz}$. The lower bound in the
``default-type boundary'' case $\aexp_{\zz} > 2$ with slower rates
than in regular regression requires a completely new strategy of proof
where not only alternatives for the regression function, but also for
the error distributions are tested against each other.

In Section~\ref{5} we provide some numerical simulations in order to
evaluate the finite sample performance of the estimator. Smaller values
in $\aexp_{\zz}$ indeed lead to significantly improved estimation
results. The bandwidth selection shows a quite different behavior from
the regular regression case due to taking local extremes. Applications
to empirical data from sunspot observations and annual best running
times on 1500~m are presented. Most proofs are deferred to Section~\ref{6}, and auxiliary lemmas and details regarding the sharpness
estimation are given in the supplementary material~\cite{jirmeireisssuppl}.

\section{Methodology} \label{2}

Our approach is a local polynomial estimation based on local extreme
value statistics. We fix some $x\in[0,1]$ and consider the
coefficients $(\hat{b}_j)_{j=0,\ldots,\beta^*}$ which minimize
the objective function
%
\begin{equation}
\label{defnestimator} (b_0,\ldots,b_{\beta^*}) \mapsto\sum
_{\llvert x_i-x\rrvert\leq h_k} \sum_{j=0}^{\beta^*}
b_j (x_i-x)^j,
\end{equation}
under the constraints $Y_i \leq\sum_{j=0}^{\beta^*} b_j
(x_i-x)^j$ for all $i$ with $\llvert x_i-x\rrvert\leq h_k$. Set
$\tilde{f}_k(x):=
\hat{b}_0$. As an estimator of $f$ we define
%
\begin{equation}
\label{eqvartheta} \tilde{f}_k\dvtx  x \mapsto\tilde{f}_k(x),
\qquad x\in[0,1],
\end{equation}
where the bandwidth $h_k>0$ remains to be selected.

\begin{figure}

\includegraphics{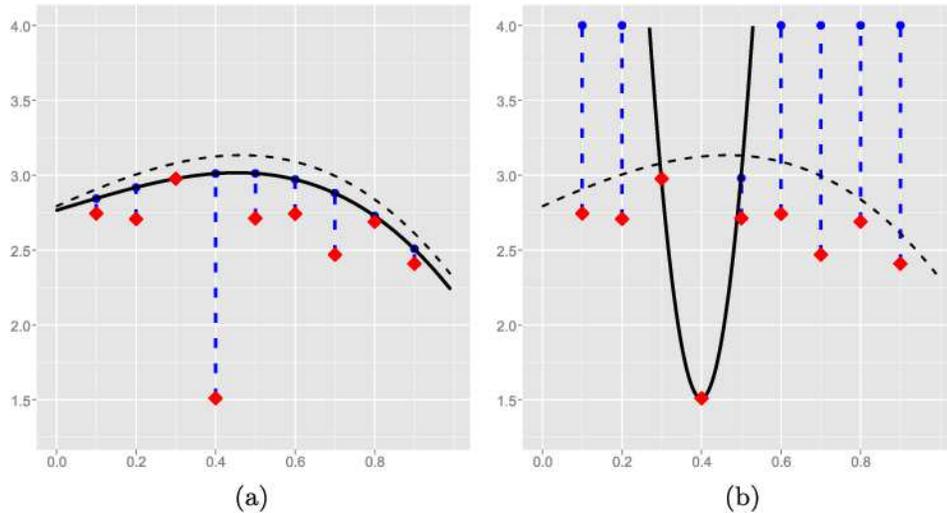}

\caption{Dashed line: true function;
solid line: estimated function;
red squares: sample points;
estimation point: fourth sample point from the left.}\label{figcompareest}\label{figunser}\label{fighall}
\end{figure}

If $-\varepsilon_j$ is exponentially distributed and the regression
function a polynomial of maximal degree $\beta^*$ on the interval
$[x-h_k,x+h_k]$, then $\tilde{f}_k(x)$ is the maximum likelihood
estimator (MLE), whence the approach can be seen as a local quasi-MLE
method; see also Knight~\cite{knight2001}. The idea of local
polynomial estimators in frontier estimation was already employed, for
instance, in Hall et al.~\cite{hall1998}, Hall and Park~\cite
{hall2004} and Hall and van Keilegom~\cite{hallkeilegom2009}.
However, in contrast to their local linear estimators (and their higher
order extensions), the sum over the evaluations at $x_i$ in the
neighborhood of $x$ is minimized instead of the area. This marks a
substantial difference and is crucial for our setup.
Already in the case of quadratic polynomials~$p$, it might occur that
the minimization of just the value $p(x)$ under the support constraints
yields the inappropriate estimator $\tilde f_k(x)=-\infty$ if $x$ is
not a design point and $\tilde f_k(x)=Y_i$ if $x=x_i$ because a sufficiently
steep parabola always fits the constraints. This problem is visualized
in Figure~\ref{figcompareest}, where Figure~\ref{figunser}(a)
corresponds to estimator (\ref{defnestimator}),
and Figure~\ref{fighall}(b) to the minimization approach employed in
the above references. Note that this problem may or may not occur in
practice, but it poses an obstacle for the mathematical analysis.
This is why we work with the base estimators defined in (\ref{defnestimator}).

The calculation of our estimator only requires basic linear
optimization, but its error analysis will be more involved. Note that
the formulation as a linear program is particularly important for
implementation purposes, since our adaptive procedure requires the
computation of many sequential estimators as the bandwidth $h_k$ increases.

%

The adaptation problem consists of finding an (asymptotically) optimal
bandwidth $h_k$ when neither the regression function $f$ nor the
specific boundary behavior of the errors $(\varepsilon_j)$ is known,
which leads to different convergence rates. We follow the method
inaugurated by Lepski~\cite{lepski1990} and consider geometrically
growing bandwidths with $h_0 = n^{\mathfrak{h}_0 - 1}$, $\mathfrak
{h}_0 \in(0,1)$ and
%
\begin{eqnarray}\label{eqbandwidth}
{h}_k &=& {h}_0 \rho^k,\qquad k = 0,\ldots,K+1
\nonumber\\[-10pt]\\[-10pt]
\eqntext{\mbox{where }\rho> 1\mbox{ and }K = \bigl\lfloor
\log_{\rho} \bigl(n^{1 - \mathfrak
{h}_0} \bigr) \bigr\rfloor.}
\end{eqnarray}
The purely data-driven estimator $\hat{f}:=\tilde
{f}_{\hat{k}}$ is defined as
%
\begin{equation}
\label{defnleskiestimator} \hat{k}:=\inf\bigl\{k = 0,\ldots,K |
\exists l \leq k\dvtx
\llVert\tilde{f}_{k+1} - \tilde{f}_l \rrVert> \hat{\z}_l^T + \hat{\z}_{k+1}^T
\bigr\} \wedge K.
\end{equation}
The critical values $\hat{\z}_l^T$, $l=0,\ldots,K+1$
depend on the observations $\{Y_i\}_{1 \leq i \leq n}$, and will be
specified below. The basic idea is to increase the bandwidth $h_k$ as
long as the distance (in some suitable seminorm $\llVert\cdot\rrVert
$) between
the estimators is not significantly larger than the usual stochastic
fluctuations of the estimators such that at $\hat k$ the bias is not
yet dominating.
In order to choose $\hat{\z}_l^T$, the extreme-value index $\aexp
_{\zz}$, $\bd_{\zz}$ and the constant $\cf_{\zz}$ from equations
(\ref{defU}) and (\ref{defL}) have to be estimated locally. For that
purpose a quasi-negative-Hill method is developed in Section~\ref{32}.


\section{Asymptotic upper bounds} \label{3}

In this section we will study the convergence rate of our estimator
$\hat{f} = \tilde{f}_{\hat{k}}$ with $\hat{k}$ as
defined in (\ref{defnleskiestimator}) when the sample size $n$ tends
to infinity. We will consider both the pointwise risk $\mathbb{E}_f
\llvert\hat{f}(x) - f(x)\rrvert^2$ for some fixed $\zz\in[0,1]$
and the $L_q$-risk
$\mathbb{E}_f \int_0^1 \llvert\hat{f}(x) - f(x)\rrvert^q \,dx$ for
$q \geq1$. To deal with the upper bounds, first some preparatory
remarks and work are necessary. Throughout this section, we suppose:
%
\begin{ass}\label{assmain}
(i) $\cf_x, \bd_x, \aexp_x \in H_{[0,1]}(\beta_0,L_0)$,
where $\beta_0, L_0 > 0$ and $\inf_{x \in[0,1]} \aexp_x,\cf_x > 0$,

(ii) $\max_{1 \leq j \leq n}\mathbb{E} [\llvert\varepsilon
_j\rrvert] < \infty$,

(iii) $(\varepsilon_j)$ are independent, and the distribution
of $\varepsilon_j$ satisfies (\ref{defU}), (\ref{defL}).
\end{ass}
For our theoretical treatment, an important quantity in the sequel is
the approximative tail-function
%
\begin{equation}
\label{defnAfisrt} A_x(y) = -\cf_x y^{-1/\aexp_x}
\log(y)^{\bd_x},
\end{equation}
since it asymptotically describes the quantile $\mathcal{U}_{x}(y)$.

\subsection{General upper bounds} \label{31}

Most of our analysis relies on Theorem~\ref{teocon} and Proposition~\ref{propestimator} below. These give rise to a decomposition, where
the error for the implicitly defined base estimators $\tilde{f}_k(x)$
in (\ref{eqvartheta}) is split\vspace*{2pt} into a deterministic and a~stochastic
error part. Even though $\tilde{f}_k(x)$ is highly nonlinear, we
obtain a relatively sharp and particularly simple upper bound.

%
\begin{teo} \label{teocon}
For any $x\in[0,1]$ and $\beta_{\zz}\in(0,\beta^*+1]$ there exist
constants $c(\beta^*,L_{\zz})$, $c(\beta^*)$ and $J(\beta^*)$, only
depending on $\beta^*$ and $L_{\zz}$, respectively, such that
\begin{eqnarray*}
&& \bigl\llvert\tilde{f}_k(x) - f(x)\bigr\rrvert
\\
&&\qquad \leq c\bigl(
\beta^*,L\bigr) h_k^{\beta_{\zz}}
\\
&&\quad\qquad{}+ c\bigl(\beta^*\bigr)\max\bigl\{ \bigl\llvert Z_j(h_k,x)
\bigr\rrvert\dvtx  j=1,\ldots,2J\bigl(\beta^*\bigr), x+h_k{\mathcal
I}_j \subseteq[0,1] \bigr\}
\end{eqnarray*}
holds true for all $f \in H_{\NN(\zz)}(\beta,L)$ where
\begin{eqnarray*}
Z_j(h_k,x) &:=& \max\{\varepsilon_i\dvtx
x_i \in x + h_k {\mathcal I}_j \}\quad\mbox{and}
\\
{\mathcal I}_j &:=& \bigl[-1+(j-1)/J\bigl(\beta^*\bigr),-1+j/J\bigl(
\beta^*\bigr)\bigr].
\end{eqnarray*}
\end{teo}

%
\begin{rem}\label{remconcentration}
Interestingly, this decomposition holds true for any underlying
distribution function $F$ and dependence structure within $(\varepsilon
_j)$. Its proof is entirely based on nonprobabilistic arguments and
has an interesting connection to algebra. A generalization to arbitrary
dimensions or other basis functions than polynomials seems challenging.
\end{rem}

We continue the range of the indices $j$ of the $(x_j,\varepsilon_j)$
from $\{1,\ldots,n\}$ to $\mathbb{Z}$ while the equidistant location
of the $x_j$ and the independence of the $\varepsilon_j$ is
maintained. Then, Theorem~\ref{teocon} yields that with $c^* = c(\beta^*,L)$
%
\begin{eqnarray}
\label{eqpointwise} \bigl\llvert\tilde{f}_k(x) - f(x)\bigr\rrvert&\leq&
c^* h_k^{\beta_{\zz}} + c\bigl(\beta^*\bigr) \max\bigl\{ \bigl
\llvert Z_j(h_k,x)\bigr\rrvert\dvtx  j=1,\ldots,2J\bigl(
\beta^*\bigr) \bigr\},\hspace*{-30pt}
\\
\label{eqL2setup} \llVert\tilde{f}_k - f\rrVert_q &\leq&
c^* h_k^{\beta_{\zz}} + c\bigl(\beta^*\bigr) \bigl\llVert\max\bigl
\{ \bigl\llvert Z_j(h_k,\cdot)\bigr\rrvert\dvtx  j=1,\ldots,2J\bigl(\beta
^*\bigr) \bigr\}\bigr\rrVert_q,\hspace*{-30pt}
\end{eqnarray}
where $\llVert\cdot\rrVert_q$ denotes the $L_q([0,1])$-norm, $q\geq
1$. To
pursue adaptivity, suppose that in terms of some seminorm $\llVert\cdot
\rrVert$,
we can bound the error via
%
\begin{equation}
\label{eqcond1} \llVert\tilde{f}_k - f\rrVert\leq{R}_k
+ B_k \qquad\forall k=0,\ldots,K+1, f\in H_{\NN(\zz)}(\beta,L),
\end{equation}
for some nonnegative random variables $B_k,R_k$, where $B_k$ increases
in $k$ and ${R}_k$ decreases in $k$. Neither the $B_k$ nor the ${R}_k$
depend on $f$, only on $\beta_{\zz}$ and $L_{\zz}$. In the sequel
$B_k$ will be a bias upper bound while ${R}_k$ is a bound on the
stochastic error, which here---in contrast to usual mean
regression---decays in $k$ for each noise realisation. The following
fundamental proposition addresses both, the pointwise and the
$L_q$-risk of the adaptive estimator, since the pointwise distance of
function values at some $x$ as well as the $L_q$-distance of functions
on $[0,1]$ define seminorms for $q \geq1$.

%
\begin{prop}\label{propestimator}
Let $\llVert\cdot\rrVert$ denote some seminorm, and let $\tilde
{f}_k$, $f$ lie
in the corresponding normed space. Assume (\ref{eqcond1}) and that the
$\hat{\z}_{k}^T$ decrease a.s. in $k$. Defining the oracle-type index
%
\begin{equation}
\label{defnoracleestimator} \hat{k}^*:= \inf\bigl\{k = 0,\ldots,K-1\dvtx  B_{k+1}
> \hat{\z}_{k+1}^T/2 \bigr\} \wedge K,
\end{equation}
we obtain for $q \geq1$:
\begin{longlist}[(a)]
\item[(a)] $\displaystyle \mathbb{E}_f \bigl[\llVert\hat{f} - \tilde
{f}_{\hat{k}^*}\rrVert^q \ind\bigl(\hat{k} > \hat{k}^*
\bigr) \bigr]^{1/q} \leq\mathbb{E}_f \bigl[\bigl(\hat{\z}_{\hat
{k}^*}^T\bigr)^q \bigr]$,\vspace*{3pt}

\item[(b)] $\displaystyle \mathbb{E}_f \bigl[\llVert\hat{f} - \tilde
{f}_{\hat{k}^*}\rrVert^q \ind\bigl(\hat{k} < \hat{k}^*
\bigr) \bigr]^{1/q}$
\\[7pt]
$\displaystyle\hspace*{30pt}\qquad \leq2^{(2q-1)/q} \mathbb{E}_f
\bigl[\hat{\z}_{\hat{k}^*}^q\bigr]^{1/q}$
\\[7pt]
$\displaystyle\hspace*{30pt}\quad\qquad{} + 2^{(2q-1)/q}\sum_{k=0}^{K-1}
\mathbb{E}_f \bigl[{R}_k^q \ind\bigl(
\exists l \leq k\dvtx  {R}_l > \hat{\z}_l^T/2
\bigr) \bigr]^{1/q}$.
\end{longlist}
\end{prop}

\subsection{Critical values and their estimation} \label{32}

Our adaptive procedure and particularly the question of optimality
crucially hinge on the (estimated) critical values~$\hat{\z}_k^T$, and
thereby as a quantile for the distribution function $F_{\zz
}(-y^{-1})$ as $y \to\infty$. In the literature (de Haan and
Ferreira~\cite{dehaanbook2006}), the standard, nonparametric quantile
estimator is constructed via the approximation
%
\begin{equation}
\label{eqquantfails} \mathcal{U}_{\zz} (t y ) \approx\mathcal{U}_{\zz}
(y ) + a_{\zz} (t ) \bigl(y^{-1/\aexp_{\zz}} - 1 \bigr)
\aexp_{\zz}\qquad\mbox{as }t,y \to\infty,
\end{equation}
where the function $a_{\zz} (t )$ is a so-called \textit
{first-order scale function}. Unfortunately, this approach fails in our
setup. The reason for this failure is the severely shifted sample
$(Y_j)$ (we do not observe $\varepsilon_j$) and the particular type of
interpolation used in~(\ref{eqquantfails}), which leads to an
insufficient rate of convergence in the above approach. The bias that
is induced by the shift will be present in any estimation method. This
fact makes us believe that under model (\ref{eq11}), quantile
estimation for general regular varying distributions is not possible.
Since for any $t > 0$ we have the relation
%
\begin{eqnarray}
\label{eqmotivateA} F_{\zz} \bigl(A_{\zz}(n/t)
\bigr)^n \bigl(1 + \mbox{\scriptsize$\mathcal{O}$}(1) \bigr) &=&
F_{\zz} \bigl(\mathcal{U}_{\zz
}(n/t) \bigr)^n\nonumber
\\
&=& (1 - t/n )^n
\\
&=& e^{-t} \bigl(1 + \mbox{\scriptsize$
\mathcal{O}$}(1) \bigr)\qquad\mbox{as }n \to\infty,\nonumber
\end{eqnarray}
a viable alternative is provided by a plug-in estimator $\widehat
{A}_{\zz}(y)$, based on suitable estimates $\hat{\aexp}_{\zz}$,
$\hat{\bd}_{\zz}$ and $\hat{\cf}_{\zz}$. Here, the shift
may be overcome by location invariant estimators for these quantities.
The fact that these parameters additionally vary in~$\zz$ with unknown
smoothness degree adds another level of complexity and needs to be
dealt with in a localized, adaptive manner. At this stage, it is worth
mentioning that our adaptive procedure does not hinge on any particular
type of quantile estimator. As a matter of fact, we only require the
following property of an admissible quantile estimator $\widehat{A}_{\zz}(y)$.
%
\begin{defn}\label{defnquantestadmissible}
Given $\zz\in[0,1]$, let $\mathcal{Y}_{\zz} = [\log n, n^{4/\aexp
_{\zz}}]$ and $s \in\{0,1\}$. We call $\widehat{A}_{\zz}(y)$
admissible if for any fixed $v\in\N$ and constants $c_1^- < 1<c_1^+$,
which may be arbitrarily close to one, we have
\[
P_f \biggl(c_1^- \leq\sup_{y \in\mathcal{Y}_{\zz}}
\biggl\llvert\frac
{\widehat{A}_{\zz}^{(s)}(y)}{A_{\zz}^{(s)}(y)}\biggr\rrvert\leq c_1^+
\biggr) = 1 -
\OO\bigl(n^{-v} \bigr),
\]
uniformly over $f\in H_{\NN(\zz)}(\beta,L)$, where $g^{(s)}(\cdot)$
denotes the $s$th derivative of a function $g(\cdot)$.
\end{defn}

%
\begin{rem}\label{remadmL2}
Admissibility for $s = 1$ is only required in case of the $L_q$-norm loss.
\end{rem}

Now we shall construct an admissible estimator under Assumption~\ref
{assmain}. Even though the class of potential estimators seems to be
quite large under Assumption~\ref{assmain}, verifying the conditions
of Definition~\ref{defnquantestadmissible} leads to quite technical
and tedious calculations. Moreover, the requirement of location
invariance rules out many prominent estimators from the literature.
Regarding the shape parameter~$\aexp_{\zz}$, this eliminates, for
instance, Hill-type estimators as possible candidates; see Alves~\cite
{fraga2002} and de Haan and Ferreira~\cite{dehaanbook2006}. Possible
alternatives are Pickand's estimator (cf. Pickand~\cite
{pickandsIII1975} and Drees~\cite{drees1995}) or the probability
weighted moment estimator by Hosking and Wallis~\cite{hoskins1987}.
These may, however, exhibit a poor performance in practice; see, for
instance, de Haan and Peng~\cite{dehaanpeng1998} for a comparison.
In~\cite{falk1995}, Falk proposed the negative Hill estimator, which,
unlike to its positive counter part, is also location invariant; see
also de Haan and Ferreira~\cite{dehaanbook2006}. Transferring this
approach to our setup, we construct estimators $\hat{\aexp}_{\zz
}$, $\hat{\bd}_{\zz}$ and $\hat{\cf}_{\zz}$ that are
location invariant, and also inherit the favorable variance property of
Hill's estimator. Based on these estimates, we can use the plug-in estimator
%
\begin{equation}
\widehat{A}_{\zz}(y) = -\hat{\cf}_{\zz}(\log
y)^{\hat{\bd}_{\zz}} y^{-1/\hat{\aexp}_{\zz}}.
\end{equation}

To construct the estimators $\hat{\aexp}_{\zz}$, $\hat{\bd}_{\zz}$ and
$\hat{\cf}_{\zz}$ for fixed $\zz\in[0,1]$,
consider the neighborhoods $\NN_k(\zz) = \{y\dvtx  \llvert\zz-y\rrvert
\leq
h_k\}$ for $k = 0,\ldots,K-1$. Introduce the sets $\mathcal
{S}_{k}(\zz) = \{Y_i\dvtx  i/n \in\NN_{k}(\zz)\}$, and note that
its cardinality $\bar n_k(\zz):= \# \mathcal{S}_{k}(\zz)$ satisfies
$n h_k \leq\bar n_k(\zz)\leq2 nh_k + 1$. Let us rearrange the sample
in $\mathcal{S}_{k}(\zz)$ as
%
\begin{equation}
Y_{1,\bar n_{k}(\zz)}, Y_{2,\bar n_{k}(\zz)},\ldots, Y_{\bar
n_{k}(\zz),\bar n_{k}(\zz)},
\end{equation}
where $Y_{j,\bar n_{k}(\zz)}$ denotes the $j$th largest $Y_{i} \in
\mathcal{S}_{k}(\zz)$. For each $k = 0,\ldots,K-1$, let $\mm_k =
\mm(\bar n_{k}(\zz) )$ such that $\mm_k/\bar n_k \to0$,
where $\bar n_k = \bar n_k(\zz)$ to lighten the notation. In the
literature, a common parametrization of $\mm_k$ is $\mm_k = {\bar
n_k}^{\mathfrak{m}}$ for $0 < \mathfrak{m}\leq1$. Before discussing
the important issue of possible choices of $\mathfrak{m}$, we formally
introduce our estimation procedure. Apart from the necessary location
invariance, an estimator of ${A}_{\zz}(y)$ should also adapt to the
unknown smoothness degree of the parameters ${\aexp}_{\zz}$, ${\bd
}_{\zz}$ and ${\cf}_{\zz}$. A related issue is dealt with in the
literature; see, for instance, Drees~\cite{drees2001} or Grama and
Spokoiny~\cite{gramaspokoiny2008}. In order to achieve this
adaptivity, we apply a Lepski-type procedure to select among
appropriate base estimators. We first tackle the problem of estimating
${\aexp}_{\zz}$. Using Falk's idea in~\cite{falk1995}, we define
%
%
\begin{equation}
\qquad\frac{1}{\hat{\aexp}_{\zz} (\mm_k )} = \frac{1}{\mm
_k}\sum_{i = 2}^{\mm_k -1}
\log\biggl(\frac{Y_{\mm_k,\bar n_k} -
Y_{1,\bar n_k}}{Y_{i,\bar n_k} - Y_{1,\bar n_k}} \biggr),\qquad k =
0,1,\ldots,K-1.
\end{equation}
Note that this estimator is clearly location invariant. For $\rho> 1$
select the index $\hat{k}_{\aexp}(\zz)$ via
%
\begin{eqnarray}\label{defnaestimatorindex}
\hat{k}_{\aexp}(\zz) &:=&\inf\bigl\{k = 0,\ldots,K -1 | \exists l \leq k\dvtx
\nonumber\\[-8pt]\\[-8pt]
&&\phantom{\inf\bigl\{}
\bigl\llvert\hat{\aexp}_{\zz}^{-1} (\mm
_{k+1} ) - \hat{\aexp}_{\zz}^{-1} (
\mm_l )\bigr\rrvert> \rho^{-k} (\log n)^{-1}
\bigr\} \wedge K.\nonumber
\end{eqnarray}
As a final estimator, we put
%
\begin{equation}
\hat{\aexp}_{\zz}^{-1} = \hat{\aexp}_{\zz}^{-1}
(\mm_{\hat{k}_{\aexp}} )\qquad\mbox{where }\hat{k}_{\aexp
} =
\hat{k}_{\aexp}(\zz).
\end{equation}

For the estimation of ${\bd}_{\zz}$, we proceed in a similar manner.
For $k = 0,1,\ldots,K-1$, we put
%
\begin{equation}
\hat{\bd}_{\zz} (\mm_k ) = \frac{1}{\mm_k \log\log
\bar n_k}\sum
_{i = 2}^{\mm_k -1}\log\biggl(\frac{Y_{i,\bar n_k} -
Y_{1,\bar n_k}}{(\bar n_k/i)^{-1/\hat{\aexp}_{\zz}(\mm_k)} -
(\bar n_k/1)^{-1/\hat{\aexp}_{\zz}(\mm_k)}}
\biggr),\hspace*{-35pt}
\end{equation}
and select the index $\hat{k}_{\bd}(\zz)$ via
%
\begin{eqnarray}
\label{defnbestimatorindex}
\hat{k}_{\bd}(\zz)&:=&\inf\bigl\{k = 0,\ldots, \hat{k}_{\aexp}(\zz) | \exists l \leq k\dvtx
\nonumber\\[-8pt]\\[-8pt]
&&\hphantom{\inf\bigl\{}
\bigl\llvert\hat{\bd}_{\zz
} (\mm_{k+1} ) - \hat{\bd}_{\zz} (\mm_l )\bigr\rrvert> \rho^{-k} (
\log\log n)^{-1} \bigr\} \wedge K.\nonumber
\end{eqnarray}
As final estimator, we then put
%
\begin{equation}
\hat{\bd}_{\zz} = \hat{\bd}_{\zz} (
\mm_{\hat{k}_{\bd}} )\qquad\mbox{where }\hat{k}_{\bd} = \hat{k}_{\bd}(\zz).
\end{equation}
Interestingly, it turns out that $\hat{\bd}_{\zz} = \bd_{\zz}
+ \frac{\log\cf_{\zz}}{\log\log n h_{\hat{k}_{\bd}}} (1
+ \mbox{\scriptsize $\mathcal{O}$}_P(1) )$. Since this implies that
\begin{eqnarray*}
(\log nh_k)^{\hat{\bd}_{\zz}} &=& \cf_{\zz}(\log
nh_k)^{\bd
_{\zz}} \bigl(1 + \mbox{\scriptsize$
\mathcal{O}$}_P(1) \bigr)\qquad\mbox{for }k = 0,\ldots, K-1,
\end{eqnarray*}
there is no need to specifically estimate $\cf_{\zz}$, it is included
in the bias for free. We are thus lead to the definition of our estimator
%
\begin{equation}
\label{defnAestimate} \widehat{A}_{\zz}(y) = -(\log y)^{\hat{\bd}_{\zz}}
y^{-1/\hat{\aexp}_{\zz}}.
\end{equation}

For the consistent estimation of $\widehat{A}_{\zz}(\cdot)$, we need
a relation between the initial bandwidth $h_0$ and the bias, induced by
the parameter $\beta_{\zz}$. Note that such an assumption is
inevitable, since any adaptive estimation procedure needs to start off
with some initial bandwidth. Thus in the sequel, we will assume that
%
\begin{equation}
\label{eqparamrelation0} h_0^{\beta_{\zz}} \bigl\llvert A_{\zz} (
\mm_0 )\bigr\rrvert^{-1} = \mbox{\scriptsize$\mathcal{O}$}
\bigl((\log n)^{-1} \bigr).
\end{equation}
If $\mathfrak{h}_0,\mathfrak{m}> 0$ is such that
%
\begin{equation}
\label{eqparamrelation1} \mathfrak{m}\mathfrak{h}_0 < (1 -
\mathfrak{h}_0)\aexp_{0} \beta_0,
\end{equation}
for some lower bounds
%
\begin{equation}
\label{eqlowbnd} \aexp_0\leq\aexp_{\zz}\quad\mbox{and}\quad
\beta_0 \leq\beta_{\zz}
\end{equation}
on the unknown parameters, then (\ref{eqparamrelation0}) is valid. In
the supplementary material~\cite{jirmeireisssuppl}
we prove the following result under the more general
Assumption 10.1, which is implied by Assumption~\ref{assmain}.
%
\begin{prop}\label{propestimateunivquantile}
Grant Assumption 10.1, and suppose that (\ref{eqparamrelation0}) is
valid. Then $\widehat{A}_{\zz}(y)$ defined in (\ref{defnAestimate})
is admissible.
\end{prop}

In practice the negative Hill estimator works well for $\aexp_{\zz}
\in(0,3/2]$ (and small~$\bd_{\zz}$), but has increasing
(asymptotically negligible) bias for $\aexp_{\zz} > 3/2$ and $\bd
_{\zz} \neq0$, which should be corrected in applications; see also
Section~\ref{5}, paragraph~(B). Also note that our assumptions in Assumption~\ref{assmain} include cases where a CLT for an estimator $\hat{\aexp
}_{\zz}$ fails to hold, and only slower rates of convergence than
$m_k^{-1/2}$ are possible. This is particularly the case if $\bd_{\zz
} \neq0$; we refer to de Haan and Ferreira~\cite{dehaanbook2006} for details.
In practice, the choice of the actual bandwidth $\mm_k$ (and hence
$\mathfrak{m}$) is of significant relevance, and much research has
been devoted to this subject; see, for instance, Drees~\cite
{Drees1998} and Drees et al.~\cite{drees2001how}. In~\cite
{fraga2002}, Alves addresses 
this question for a related (positive) location invariant Hill-type
estimator both in theory (Theorem~2.2) and practice (concluding remarks
and algorithm). Transferring the practical aspects, this amounts to the
choice $\mm_k = 2{\bar n_k}^{\mathfrak{m}}$, $\mathfrak{m}= 2/3$ in
our case. Still, any other choice also leads to the total optimal rates
presented in Theorems~\ref{Tpointwise} and~\ref{TL2upperbound}, as
long as $0 < \mathfrak{m}< 1$ holds.

\subsection{Pointwise adaptation} \label{secpointwise}

Throughout this subsection we fix a point $\zz\in[0,1]$. For the
seminorm in Proposition~\ref{propestimator} we take $\llVert f\rrVert:=\llvert f(\zz
)\rrvert$. According to Theorem~\ref{teocon}, we set
%
\begin{eqnarray}\label{EqRk2}
\nonumber
B_k &:=& c\bigl(\beta^*,L\bigr) h_k^\beta,
\nonumber\\[-8pt]\\[-8pt]\nonumber
{R}_k &:=& c\bigl(\beta^*\bigr) \max\bigl\{ \bigl\llvert
Z_j(h_k,\zz)\bigr\rrvert\dvtx  j=1,\ldots,2J\bigl(\beta^*
\bigr) \bigr\},
\end{eqnarray}
in the notation of (\ref{eqcond1}). The nonnegativity and monotonicity
constraints on $B_k$ and $R_k$ are satisfied since $h_k$ increases. We
define the oracle and estimated critical values as
\begin{eqnarray}
\nonumber
\z_k(\zz) &=& 4 c\bigl(\beta^*\bigr)\biggl\llvert
A_{\zz} \biggl(\frac{\aexp_{\zz}
n h_k}{4 J(\beta^*) \log n} \biggr)\biggr\rrvert,\qquad\hat{\z}_k(\zz)
= 4 c\bigl(\beta^*\bigr)\biggl\llvert
\widehat{A}_{\zz} \biggl(\frac
{\hat{\aexp}_{\zz}n h_k}{4 J(\beta^*)\log n } \biggr)\biggr\rrvert,
\end{eqnarray}
for $k=0,\ldots,K-1$ and set ${\z}_K(\zz) = \hat{\z}_K(\zz):=0$. To
lighten the notation, we often drop the index $\zz$ and write
$\z_k$ and $\hat{\z}_k$. As outlined earlier in (\ref
{eqmotivateA}), this definition is motivated by the fact that $\mathcal
{U}_{\zz}(y) \approx A_{\zz}(y)$ as $y \to\infty$. The critical
values can thus be viewed as an appropriate estimate for certain
extremal quantiles. The additional $\log n$-factor turns out to be the
price to pay for adaption. We proceed by introducing the estimated
truncated critical values as
%
\begin{equation}
\label{eqSchaetzerzeta} \hat{\z}_k^T = \min\{\hat{\z}_k,1 \}.
\end{equation}
The truncation of the estimator $\hat{\z}_k^T$ is required to
exclude a possible pathological behavior both in theory and practice.
Note that this does not affect its proximity to $\z_k$ if $\hat{\z}_k$
is consistent, since $\z_k \to0$ uniformly in $k = 0, \ldots,K-1$. We have the following pointwise result.


%
\begin{teo}\label{Tpointwise}
Fix $\zz\in[0,1]$, and suppose $\aexp_{\zz}, \bd_{\zz}, \cf_{\zz
}$ and $\beta_{\zz}\in(0,\beta^*+1]$ are unknown with $\mathfrak
{h}_0 < \beta_{\zz} \aexp_{\zz} / (\beta_{\zz} \aexp_{\zz} +
1)$. If Assumption~\ref{assmain} holds, then
\begin{eqnarray*}
&& \sup_{f \in H_{\NN(\zz)}(\beta,L)} \mathbb{E}_f \bigl[ \bigl(
\hat{f}(\zz) - f(\zz) \bigr)^2 \bigr]
\\
&&\qquad = \OO\bigl((n /\log
n)^{(-2 \beta_{\zz})/(\aexp_{\zz} \beta_{\zz}+1)} (\log n)^{(2\aexp_{\zz} \bd_{\zz} \beta_{\zz})/(\aexp_{\zz}
\beta_{\zz}+1)} \bigr).
\end{eqnarray*}
\end{teo}

As will be demonstrated in Section~\ref{seclowerbounds}, this result
is optimal in the minimax sense.

\subsection{$L_q$-adaptation} \label{secL2}

Let us consider the $L_q([0,1])$-norm as seminorm in Proposition~\ref
{propestimator}. Due to (\ref{eqL2setup}) we can choose
%
\begin{eqnarray}
B_k &:=& c\bigl(\beta^*,L\bigr) h_k^\beta,
\label{EqBk}
\\
R_k &:=& c\bigl(\beta^*\bigr) \bigl\llVert\max\bigl\{ \bigl\llvert
Z_j(h_k,\cdot)\bigr\rrvert\dvtx  j=1,\ldots,2J\bigl(\beta^*
\bigr) \bigr\}\bigr\rrVert_q \label{EqRk}
\end{eqnarray}
in the notation of (\ref{eqcond1}). We verify that the nonnegativity
and monotonicity constraints on $B_k$ and $R_k$ are satisfied for $\rho
>1$ in (\ref{eqbandwidth}) since for any $x\in[0,1]$ each interval $x
+ h_k {\mathcal I}_j$, $j=1,\ldots,2J(\beta^*)$ is included in $x +
h_{k+1} {\mathcal I}_{j'}$ for some $j'=1,\ldots,2J(\beta^*)$ for any
$k$. Throughout this paragraph, we assume that the parameters $\aexp
_{x}, \bd_{x}, \cf_{x}$ remain constant for $x \in[0,1]$. We denote
these with $\aexp_{F}, \bd_{F}, \cf_{F}$, and the corresponding
$A_x(\cdot)$ with $A_F(\cdot)$.

The construction of the critical values is more intricate compared to
the pointwise case, and relies on the following quantity. Introduce
%
\begin{equation}
\qquad \widehat{\mathit{IU}}_n(s,q)= \biggl(\int_{n^{-2/\hat{\aexp
}_F}}^{n^{1/2}}
\bigl((-\widehat{A}_F)^q (s/y ) \bigr)^{(1)}
\exp(-y) \,dy \biggr)^{1/q}, \qquad q \geq1,
\end{equation}
and the corresponding version ${\mathit{IU}}_n(s,q)$ where we replace $\widehat
{A}_F (\cdot)$ by ${A}_F (\cdot)$ and
$\hat{\aexp}_F$ by $\aexp_F$ [recall that $g^{(s)}(\cdot)$ denotes the
$s$th derivative of a function $g(\cdot)$]. For $k=0,\ldots,K-1$, we
introduce the critical values as
\begin{eqnarray}
\nonumber
\z_k &=& \sqrt{5} c\bigl(\beta^*\bigr)\biggl\llvert
{\mathit{IU}}_n \biggl(\frac{n h_k}{6
J(\beta^*)},q \biggr)\biggr\rrvert,\qquad
\hat{\z}_k = \sqrt{5} c\bigl(\beta^*\bigr)\biggl\llvert
\widehat{\mathit{IU}}_n \biggl(\frac{n h_k}{6 J(\beta
^*)},q \biggr)\biggr\rrvert,
\end{eqnarray}
and set ${\z}_K = \hat{\z}_K:=0$. Moreover, we define the
corresponding truncated values as
%
\begin{equation}
\hat{\z}_k^T = \min\{\hat{\z}_k,1
\}.
\end{equation}
Unlike the pointwise case, the critical values do not correspond to an
extremal quantile, but they can be considered as an estimate of
$\mathbb{E}[R_k^q]^{1/q}$. This already indicates that the $L_q$-case
is substantially different from the pointwise situation, and indeed
additional, more refined arguments are necessary to prove the result
given below.
%
\begin{teo} \label{TL2upperbound}
Suppose $\aexp_F>0$ and $\beta\in(0,\beta^*+1]$ are unknown with
$\beta\aexp_F / (\beta\aexp_F + 1)>\mathfrak{h}_0$. We select
$\rho\in\N$ with $\rho>1$. If $q \geq1$, then the adaptive
estimator $\hat{f}$ from Section~\ref{2} satisfies
\[
\sup_{f\in H_{[0,1]}(\beta,L)} \mathbb{E}_f \bigl[\llVert
\hat{f} - f\rrVert_q^q \bigr] = \OO
\bigl(n^{(-q\beta)/(\aexp_{F} \beta_+1)} (\log n)^{(q \beta\aexp_{F} \bd_{F})/(\aexp_{F} \beta
_+1)} \bigr). %
\]
\end{teo}

%
\begin{rem}\label{remL2case}
If one allows for $\aexp_{x}, \bd_{x}, \cf_{x} \in H (\beta
_0,L )$ for $x \in[0,1]$, the above result remains valid if one
takes the supremum over the above bound. This result is also optimal in
the minimax sense.
\end{rem}

Theorem~\ref{TL2upperbound} shows that the estimator $\hat{f}$ is
$L_q$-adaptive; that is, it attains the minimax rates, which are
optimal in the oracle setting of known $\aexp_F$ and $\beta$,
although it does not use these constants in its construction; see
Theorem~\ref{T421} below for the lower bound.

\section{Asymptotic lower bounds} \label{seclowerbounds}

We show that the logarithmic loss in the convergence rate in Theorem
\ref{Tpointwise} is unavoidable with respect to any estimator sequence
of $f$. First, we treat the case $\aexp_x\in(0,2]$ for which we
derive a lower bound, even for a known error distribution. It suffices
to treat the case where $\aexp= \aexp_x$, $\bd= \bd_x$ and $\cf=
\cf_x$ remain constant for $x \in[0,1]$. We maintain this convention
throughout this section.

We assume that the $\varepsilon_j$ have a Lebesgue density
$f_\varepsilon$ which is continuous and strictly positive on $(-\infty,0)$,
and vanishes\vspace*{1pt} on $[0,\infty]$. Moreover, we impose that the $\chi^2$-distance for the parametric location problem satisfies
%
\begin{equation}
\label{eqeps} \qquad\int_{-\infty}^0 \bigl\llvert
f_\varepsilon(x+\vartheta) - f_\varepsilon(x)\bigr\rrvert^2
/ f_\varepsilon(x) \,dx \leq c_\varepsilon\vartheta^{\aexp} \llvert
\log\vartheta\rrvert^{-\aexp\bd} \qquad\forall\vartheta\in(0,1),
\end{equation}
for some $\aexp\in(0,2]$ and $\bd\in\mathbb{R}$. Note
that $\aexp$ and $\bd$ correspond to $\aexp_x$ and $\bd_x$,
respectively, in (\ref{defU}) and (\ref{defL}) with uniform $x$. As
examples for such error densities with $\bd=0$, we consider the
reflected gamma-densities
\[
f_{\lambda}(x):= \frac{1}{\Gamma(\lambda)} (-x)^{\lambda-1} \exp(x){
\mathbf1}_{(-\infty,0)}(x), \qquad x\in\mathbb{R}, %
\]
for $\lambda\in[1,2)$. Thus, by $\llvert(1+\vartheta/x)^{\lambda
-1}-1\rrvert\le\vartheta/\llvert x\rrvert$ for $x\le-\vartheta$ we have
\begin{eqnarray*}
&&\int_{-\infty}^0 \bigl\llvert f_\lambda(x+
\vartheta) - f_\lambda(x)\bigr\rrvert^2 /
f_\lambda(x) \,dx
\\
&&\qquad = \int_{-\infty}^{-\vartheta} \biggl\llvert\biggl(1 +
\frac{\vartheta}{x} \biggr)^{\lambda-1} \exp(\vartheta) - 1 \biggr
\rrvert
^2 f_\lambda(x) \,dx + \int_{-\vartheta}^0
f_\lambda(x) \,dx
\\
&&\qquad \leq2 \bigl(\exp(\vartheta)-1 \bigr)^2 + 2 \exp(2\vartheta)
\vartheta^2
\\
&&\quad\qquad{}+ 2 \exp(2\vartheta) \vartheta^2 \int
_{-1}^{-\vartheta
} x^{-2} f_\lambda(x) \,dx
+ \int_{-\vartheta}^0 f_\lambda(x) \,dx
\\
&&\qquad \leq\OO\bigl(\vartheta^2\bigr) + \frac{2}{(2-\lambda) \Gamma(\lambda
)} \exp(2
\vartheta) \vartheta^\lambda\bigl(1 - \vartheta^{2-\lambda}\bigr) +
\vartheta^\lambda/\Gamma(\lambda+1)
\\
&&\qquad = \OO\bigl(\vartheta^{\lambda}\bigr).
\end{eqnarray*}
%
Therefore, the reflected gamma-density satisfies (\ref
{eqeps}) when putting $\aexp= \lambda$. Note that (\ref{eqeps})
implies (\ref{defU}), (\ref{defL}) under the Assumption~\ref{assmain}(i). The following theorem together with the upper bound in
Theorem~\ref{Tpointwise} shows that pointwise adaptation causes a
logarithmic loss in the convergence rates, which is known from regular
regression when inserting $\aexp=2$.

%
\begin{teo} \label{Tlowerbound}
Assume condition (\ref{eqeps}), and fix some arbitrary $x_0\in[0,1]$,
$\beta_1>\beta_2>0$ and $C_0,C_1>0$. Let $\{\hat{f}_n(x_0)\}_n$
be any sequence of estimators of $f(x_0)$ based on the data $Y_1,\ldots,Y_n$ which satisfies
\[
\sup_{f\in{\mathcal H}_{[0,1]}(\beta_1,C_0)} \mathbb{E}_f \bigl\llvert
\hat{f}_n(x_0) - f(x_0)\bigr\rrvert
^2 = \OO\bigl(n^{-2\beta
_2/(1+\beta_2 \aexp)} n^{-\xi} \bigr), %
\]
for some $\xi>0$. Then this estimator sequence suffers from
the lower bound
\begin{eqnarray*}
&& \liminf_{n\to\infty} (n/\log n)^{(2\beta_2)/(1+ \aexp\beta
_2)} (\log
n)^{(-2\aexp\bd\beta_2)/(1 + \aexp\beta_2)}
\\
&&\qquad{}\times  \sup_{f\in{\mathcal
H}_{[0,1]}(\beta_2,C_1)} \mathbb{E}_f \bigl
\llvert\hat{f}_n(x_0) - f(x_0)\bigr
\rrvert^2 > 0.
\end{eqnarray*}
\end{teo}

For completeness we also derive the $L_q$-minimax optimality of the
convergence rates established by our estimator $\hat{f}$ in
Theorem~\ref{TL2upperbound}. This rectifies a conjecture after Theorem~3 in Hall and van Keilegom~\cite{hallkeilegom2009} for general
smoothness degrees.

%
\begin{teo} \label{T421}
Assume condition (\ref{eqeps}), and let $\{\hat{f}_n\}_n$ be any
sequence of estimators of $f$ based on the data $Y_1,\ldots,Y_n$.
Then, for any fixed $q\geq1$, we have
\[
\liminf_{n\to\infty} n^{\beta_2/(1+\aexp\beta_2)} (\log n)^{(-\aexp\bd\beta_2 )/(1 + \aexp\beta_2)} \sup
_{f\in
{\mathcal H}_{[0,1]}(\beta_2,C_1)} \mathbb{E}_f \bigl[\llVert\hat{f}_n
- f\rrVert_q \bigr] > 0. %
\]
\end{teo}

Now we focus on the case $\aexp>2$. To simplify some of the technical
arguments in the proofs, we restrict to the case $\bd= 0$. If $\aexp
>2$, the convergence rates become slower than in the Gaussian case.
Instead of the convenient conditions (\ref{defU})~and~(\ref{defL}),
we choose the slightly different Definition~\ref{assmainlowerbound},
under which the upper bound proofs obviously still hold true.

%
\begin{defn}\label{assmainlowerbound}
Let $\aexp>2$, $0 < \mathfrak{h}_0 \leq1$, and denote with
${\mathcal D}_n(\aexp, \mathfrak{h}_0)$ the set of all error
distribution functions whose quantile functions $\mathcal{U}^{(n)}$ satisfy:
\begin{longlist}[(ii)]
\item[(i)] $\displaystyle \sup_{y \in(0, \infty]} \biggl\llvert\frac{\mathcal
{U}^{(n)}(y)}{A(y/2)}\biggr
\rrvert\leq1\qquad\mbox{where }A(y) = - y^{-1/\aexp}$,\vspace*{9pt}

\item[(ii)] $\displaystyle \sup_n \sup_{y \in[\log N, N]} \biggl
\llvert\frac{\mathcal {U}^{(n)}(y)}{A(y)} - 1 \biggr\rrvert\llvert\log y\rrvert\leq(\log
n)^{-2}, \qquad N = n^{\mathfrak{h}_0}$.
\end{longlist}
Note that we have ${\mathcal D}_n(\aexp, \mathfrak{h}_0) \subseteq
{\mathcal D}_n(\aexp, \mathfrak{h}_0')$ if $\mathfrak{h}_0 >
\mathfrak{h}_0'$.
\end{defn}

The above conditions particularly imply that the distribution function
$F(y) = F^{(n)}(y)$ [or likewise $\mathcal{U}(y) = \mathcal
{U}^{(n)}(y)$] of the errors $\varepsilon_j$ may depend on $n$. While
the lower bound results for $\aexp\leq2$ still hold true if the error
distribution is known and independent of the design point, here two
competing types of regression errors have to be considered in the
proof. Note that the probability measure thus depends on both the
regression function $f$ and the distribution function $F$, which we
mark as $P_{f, F}$.

%
\begin{teo} \label{teo-aexp2}
Fix\vspace*{1pt} some arbitrary $x_0\in[0,1]$, $\beta_1>\beta_2> 0$ and
$C_0,C_1>0$. Let $\aexp> 2$, and suppose that $\mathfrak{h}_0 < \frac
{\beta_2}{\aexp\beta_2 + 1}$. Let $\{\hat{f}_n(x_0)\}_n$ be any
sequence of estimators of $f(x_0)$ based on the data $Y_1,\ldots,Y_n$
which satisfies
\[
\sup_{f\in{\mathcal H}_{[0,1]}(\beta_1,C_0)} \sup_{F \in{\mathcal
D}_n(\aexp, \mathfrak{h}_0)}
\mathbb{E}_{f, F} \bigl\llvert\hat{f}_n(x_0)
- f(x_0)\bigr\rrvert^2 = \OO\bigl(n^{-2\beta_2/(1+ \aexp
\beta_2)}
n^{-\xi} \bigr), %
\]
for some $\xi>0$. Then this estimator sequence suffers from
the lower bound
\[
\liminf_{n\to\infty} (n/\log n)^{(2\beta_2)/(1+\aexp\beta_2
)} \sup
_{f\in{\mathcal H}_{[0,1]}(\beta_2,C_1)} \sup_{F \in
{\mathcal D}_n(\aexp, \mathfrak{h}_0)} \mathbb{E}_{f, F}
\bigl\llvert\hat{f}_n(x_0) - f(x_0)\bigr
\rrvert^2 > 0. %
\]
\end{teo}

The proofs of Theorems~\ref{Tlowerbound},~\ref{T421} and~\ref
{teo-aexp2} are given in the supplementary material~\cite{jirmeireisssuppl}.
Theorem~\ref{T421} can be extended in a similar way to $\aexp>2$, a
detailed proof is omitted.

\section{Numerical simulations and real data application} \label{5}

The aim of this section is to highlight some of the theoretical
findings with numerical examples. We will briefly touch on the
following points:

\begin{longlist}[(A)]
\item[(A)] Performance of the estimator on different function types
and the corresponding effect on adaptive bandwidth selection.
\item[(B)] The effect of different parameters $\aexp_{\zz}$, $\bd
_{\zz}$ and $\cf_{\zz}$.
\item[(C)] Application: wolf sunspot-number.
\item[(D)] Application: yearly best men's outdoor 1500 times.
\end{longlist}

%
\begin{figure}[b]

\includegraphics{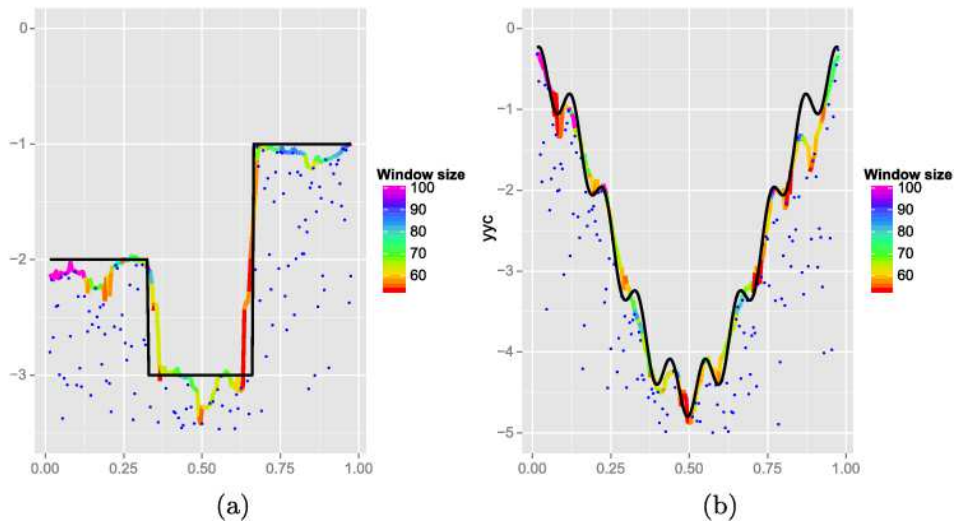}

\caption{\textup{(a)}~Function $f = f_1$, $\beta^* = 2$, $n = 200$,
$\varepsilon_j \sim \Exp(1)$;
\textup{(b)}~function $f = f_2$, $\beta^* = 2$, $n = 200$, $\varepsilon_j
\sim \Exp(1)$.}\label{figcomparemultid3}
\end{figure}

In order to illustrate the behavior of the estimation procedure, we
consider three different regression functions, displayed in black in
Figures~\ref{figcomparemultid3} and~\ref{figqq-a}(a),
\begin{eqnarray*}
f_1(x) &=& -2\cdot\ind(x < 1/3) - 3\cdot\ind(1/3\leq x < 2/3) - \ind
(2/3 < x),\qquad x\in[0,1],
\\
f_2(x) &=& -2 + 2\cos(2\pi x) + 0.3\sin(19\pi x),\qquad x\in[0,1].
\end{eqnarray*}
They are similar to those discussed in Chichignoud~\cite{chichi2012}.
Comments on the implementation and setup are given in the supplementary material~\cite{jirmeireisssuppl},
together with a\vadjust{\goodbreak} numerical comparison to oracle estimators and
additional simulations. All of the results can be reproduced by \mbox{R-code},
available at~\cite{jirmeireissrcode}.

\begin{longlist}[(A)]
\item[(A)]
Figure~\ref{figcomparemultid3} gives a first impression on the
behavior and accuracy of our estimation procedure. In both cases, the
errors $\varepsilon_j$ follow an exponential distribution $\Exp(1)$,
and the sample size is $n = 200$. The window size in Figure~\ref
{figcomparemultid3} corresponds to the local sample size, chosen by the
adaptive procedure. Even though $n$ is only of moderate size, the
estimation procedure achieves good results by essentially recovering
the shape of the underlying regression, also in the wiggly case of
function~$f_2$. Simulations of other nonparametric (adaptive)
estimators that do not take the nonregularity into account (cf. Lepski
and Spokoiny~\cite{lepskispokoiny1997} and the R-packages crs, gam,
smooth-spline, etc.) often fail to do so (with mean correction).

%
\begin{figure}

\includegraphics{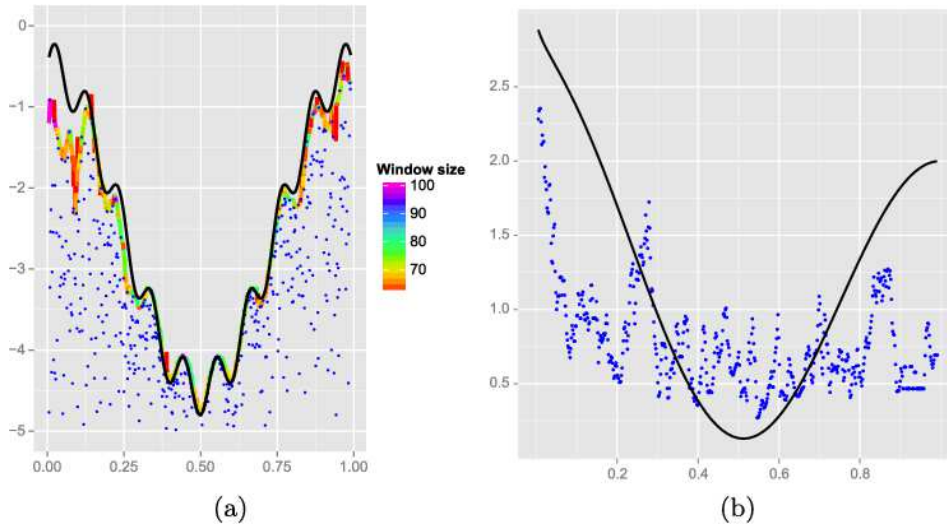}

\caption{\textup{(a)}~Function $f = f_2$, $\beta^* = 4$, $n = 600$,
$\varepsilon_j \sim\Gamma(\aexp_{\zz},1)$;
\textup{(b)}~$\aexp_{\zz}$ (black line) and $\hat{\aexp}_{\zz}$ (blue points).}\label{figcomparead3}\label{figqq-a}
\end{figure}

The effect of the shape (type) of the function on the bandwidth
selection is highlighted by a color-scheme, ranging from dark red (low)
to dark violet (high). In order to understand the ``coloring of the
estimator,'' one has to recall that the estimation procedure always
tries to fit a local polynomial which ``stays above the observations.''
At first sight, this can lead to a surprisingly large bandwidth
selection at particular spots. The bandwidth size is not necessarily an
indicator for estimation accuracy. The reason for this effect is the
maximum function: additional observations are taken into account as
long as this does not substantially change the maximum, which can lead
to a surprisingly large bandwidth selection.
\end{longlist}
\begin{longlist}[(B)]
\item[(B)]
Here, the setup is different from paragraph~(A). We consider a sample
size of $n = 600$, and we let the parameters vary in $\zz$. The impact
of $\bd_{\zz}$, $\cf_{\zz}$ (and their estimates) is rather
insignificant on the total estimator. This is not unexpected and can be
explained by the very definition in (\ref{defnAfisrt}). We therefore
focus only on the parameter $\aexp_{\zz}$ in this paragraph. We
consider the setup where the errors follow a Gamma distribution
$\varepsilon_j \sim\Gamma(\aexp_{\zz},1)$ and $\aexp_{\zz}$
varies according to the function $\aexp_{\zz} = \sin(2 \pi\zz+ \pi
/2) -\sqrt(1-(\zz-1)^2)+2$. We only discuss function $f_2$ here; a~more comprehensive comparison including an additional function $f_3$ is
given in the supplementary material~\cite{jirmeireisssuppl}. As can be clearly seen in Figure~\ref
{figcomparead3}, there is a considerable increase in estimation
accuracy as $\aexp_{\zz}$ gets closer to zero. Generally speaking,
for larger $\aexp_{\zz}$ the bias can be pronounced, and this is
indeed the case at the left top of Figure~\ref{figqq-a}(a). It simply
turns out that there are no observations at all near the regression
function $f_2$, which leads to the large gap. An approximate bias
correction [e.g., by $\widehat{\mathit{IU}}(nh_k,1)$, Section~\ref{secL2}] could be applied, but we do not pursue this any further.
Figure~\ref{figcomparead3} also reveals that the estimator $\hat{\aexp
}_{\zz}$ (blue) has a large variance and problems with quickly
oscillating regression functions $f$ (compare with the supplementary material~\cite{jirmeireisssuppl}).
On the other hand, it seems to capture the general trend of decrease
and increase to some extent. We would like to point out, however, that
these estimations are very sample dependent, and due to the relatively
small, local sample size, the actual behavior of local samples may
deviate significantly from a large sample of $\Gamma(\aexp_{\zz
},1)$-distributed random variables.
Significant overestimation leads to critical values that are too large,
which in turn results in a slight overestimation of the regression
function; see the very center of Figure~\ref{figqq-a}(a) ($0.4$ to
$0.6$). The opposite effect can be observed at both endpoints of
Figure~\ref{figqq-a}(a), where an underestimation is present, which
leads to critical values that are too small. Also note that the
negative Hill estimator generally tends to underestimate $\aexp_{\zz
}$ if $\aexp_{\zz} \geq3/2$, which is due to an (asymptotic
negligible) bias; cf. de Haan and Ferreira~\cite{dehaanbook2006}. A
thorough bias correction requires a precise second order asymptotic
expansion of the limit distribution of the negative Hill estimator,
which is beyond the scope of this paper. Note, however, that a
rudimentary bias correction is available in our implemented code.
Another, more practical option would be to consider the estimation of
$\aexp_{\zz}$ itself as a regression problem with one-sided errors,
treating the local estimates $\hat{\aexp}_{\zz}$ as ``sample.''

A similar behavior appears when considering function $f_1$, but, as can
be expected, the estimates $\hat{\aexp}_{\zz}$ are more accurate.
\end{longlist}
\begin{longlist}[(C)]
\item[(C)]
The Wolf sunspot number (often also referred to as \emph{Z\"{u}rich
number}), is a measure for the number of sunspots and groups of
sunspots present on the surface of the sun. Initiated by Rudolf Wolf in
1848 in Z\"{u}rich, this famous time series has been studied for
decades by physicists, astronomers and statisticians. The relative
sunspot number $R_t$ is computed via the formula
%
\begin{equation}
\label{eqsunspot} R_t = K_t (10 g_t +
s_t ),
\end{equation}
where $s_t$ is the number of individual spots observed at time $t$,
$g_t$ is the number of groups observed at time $t$ and $K_t$ is the
\textit{observatory factor} or \textit{personal reduction
coefficient}. The factor $K_t$ (always positive and usually smaller
than one) depends on the individual observatories around the world and
is intended to convert the data to Wolf's original scale, but also to
correct for seeing conditions and other diversions. In general, we have
the relationship
%
\begin{equation}
\label{eqmodelgen}
\mbox{observed data} = \mbox{observed fraction} \times \mbox{true value},
\end{equation}
where we always have that the random variable $\mbox{observed fraction} \in
(0,1]$. Therefore, the factor $K_t$ can be viewed as an aggregated
individual estimate for the right scaling. Over the last century, many
different models have been fit to the sunspot data; we refer to He~\cite
{heli2001} and Solanki et al.~\cite{Solanki20041084} for an overview.
In particular, the study of the sunspots has attracted people long
before 1848. Recorded observations are, for instance, due to Thomas
Harriot, Johannes and David Fabricius (in the 17th century), Edward
Maunder and many more. 
However, much uncertainty lies in these data, and the sunspot time
series before 1850 is usually referred to as ``unreliable'' or
``poor.'' It is therefore interesting to reconstruct the ``true time
series'' or at least reduce some uncertainty. We attempt do so for the
period from 1749 to 1810, based on monthly observations. Let us
reconsider model (\ref{eqsunspot}). Given~$R_t$, we may then postulate
the model
%
\begin{equation}
\label{eqsunspotpost} R_t = X_t \mathcal{S} \bigl(10
g_t^{\circ} + s_t^{\circ} \bigr),
\end{equation}
where $g_t^{\circ}$, $s_t^{\circ}$ denote the corresponding true
sunspot values, and $X_t \in(0,1]$. This means we concentrate all
random components in $X_t$, which is in spirit of model~(\ref{eqmodelgen}). We point out that this is only one possible way from a
modeling perspective; we refer to Kneip et al.~\cite{kneip2012} or
Koenker et al.~\cite{koenker1994} and the references therein for
alternatives and more general models. In our setup, the parameter
$\mathcal{S}>0$ reflects the support of the ``misjudgment'' of the
observer. For example, $\mathcal{S} \leq1$ is equivalent with the
assumption that every observer always reports less than the true value.
As we see below, it incorporates the systematic bias of the observers.
By using a $\log$-transformation, we have the additive model
%
\begin{equation}
\label{eqsunspotlog} \log R_t = \log X_t + \log
\bigl(10g_t^{\circ} + s_t^{\circ} \bigr)
+ \log\mathcal{S},
\end{equation}
which can be interpreted as a nonparametric regression problem with
stochastic error $\log X_t \in(-\infty,0]$. The goal is to estimate
the function $f(t) = \log(10g_t^{\circ} + s_t^{\circ} )$,
the ``true'' relative sunspot number. Such estimation results can serve
as input to structural physical models for sunspot activity like the
time series approaches mentioned above. Unfortunately, one can only
estimate $f(t) + \log\mathcal{S}$, where the bias $\log\mathcal{S}$
cannot be removed without any further assumptions. This is clear from
the nonidentifiability in model (\ref{eqsunspotpost}). Generally
$\mathcal{S}$ is a systematic (intrinsic) bias, which has to be
overcome using other sources of information (expert judgement). Any
other statistical approach will also suffer from such a global bias.

%
\begin{figure}[t]

\includegraphics{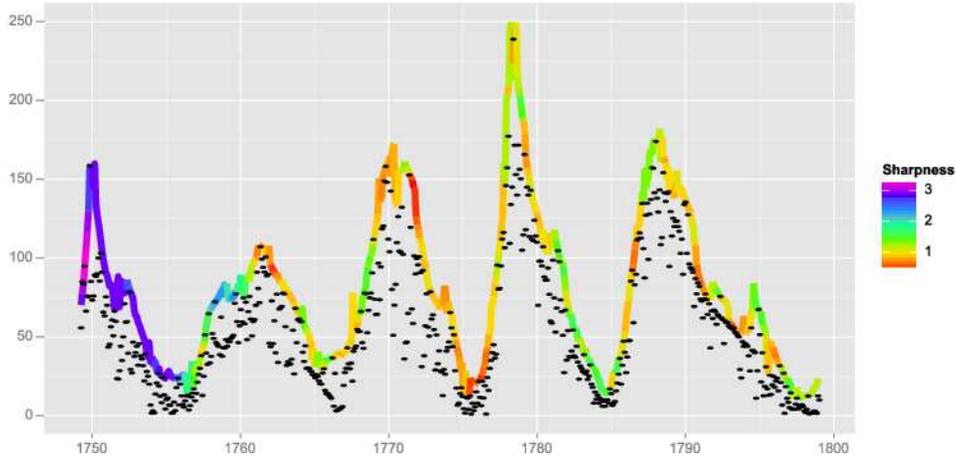}

\caption{Estimated Wolf number with $\mathcal{S} = 1$, $\beta^* =
5$.}\label{figsunspot1}
\end{figure}
%
%
\begin{figure}[b]

\includegraphics{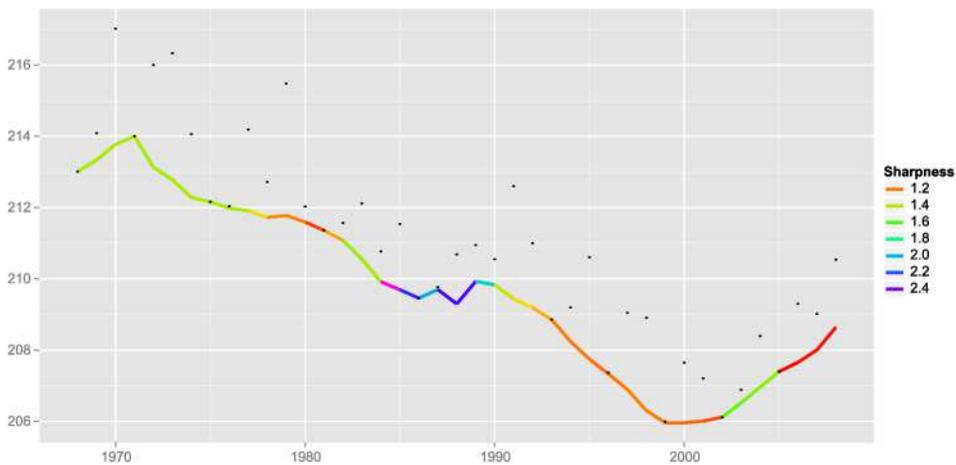}

\caption{Yearly best men's outdoor 1500~m times in seconds with
estimated boundary ($\beta^* = 2$).}\label{fig15001}
\end{figure}

The results of the estimated sunspot number is given in Figure~\ref
{figsunspot1}, where we plotted an estimate corresponding to $\mathcal
{S} = 1$. Given that observation techniques where much less advanced
and coordinated in the 18th and 19th centuries, it is\vadjust{\goodbreak} reasonable to
assume $\mathcal{S} \leq1$. Apart from the estimated sunspot number
itself, our estimation procedure provides a map from the uncertainty
level $\mathcal{S}$ to the true sunspot number~$f(t)$. The sharpness
$\aexp_{\zz}$ seems to mainly vary within the interval $[0,3.5]$.
Finally, we would like to comment on the ``peaks'' around $1768$ and
$1774$. These peaks are artifacts and originate from a too large
initial bandwith selection at these particular points. However, for the
sake of reproducibility, we have kept them and did not make any ad-hoc,
data-dependent changes.
\end{longlist}

\begin{longlist}[(D)]
\item[(D)]
As another example we discuss the yearly best men's outdoor 1500~m times
starting from 1966, depicted in Figure~\ref{fig15001} with estimated
lower boundary. Following Knight~\cite{knight2006}, the boundary can
be interpreted as the best possible time for a given year. This data
set displays an interesting behavior. As can be clearly seen from
Figure~\ref{fig15001}, the boundary steadily decreases from 1970 until
around the year 2000, followed by a sudden and sharp increase. This
event leaves room for speculation. Let us mention that until the year
2000, it had been very difficult to distinguish between the biological
and synthetical EPO. The breakthrough was achieved by Lasne and
Ceaurriz~\cite{lasneceaurriz2000}, and since then, more and more
refined and efficient doping tests have been developed. It seems
plausible that this change and advance in doping controls has lead to
the sudden increase, but it might as well be attributed to some other reason.
\end{longlist}

\section{Proof of the main results} \label{6}

Throughout the proofs, we make the following convention. For two sequences
of positive numbers $a_n$ and $b_n$, we write $a_n \gtrsim b_n$ when
$a_n \geq C b_n$ for
some absolute constant $C> 0$, and $a_n \lesssim b_n$ when $b_n \gtrsim
a_n$. Finally, we write
$a_n \thicksim b_n$ when both $a_n \lesssim b_n$ and $a_n \gtrsim b_n$ hold.

\begin{pf*}{Proof of Theorem~\ref{teocon}}
Throughout the proof, we fix some arbitrary $x\in[0,1]$ and write
$\beta= \beta_{x}$ to lighten the notation. The data $Y_i$,
$i=1,\ldots,n$, can be written as
\[
Y_i = \sum_{j=0}^{\beta^*}
b_j (x_i-x)^j + \varepsilon_i
+ \Delta_i, %
\]
where $\Delta_i:= f(x_i) - \sum_{j=0}^{\beta^*} b_j (x_i-x)^j$.
Putting $\Delta:= \max\{\llvert\Delta_i\rrvert\dvtx  \llvert x_i-x\rrvert
\leq h_k\}$, the
coefficients $b_j$ are chosen as the Taylor coefficients $b_j =
f^{(j)}(x)/j!$ for $j\leq\langle\beta\rangle$ and $b_j=0$
otherwise, such that by the H\"older condition on $f^{(\langle\beta
\rangle)}$ in the Taylor remainder term
%
\begin{equation}
\label{eqTaylor} \Delta\leq L h_k^\beta/ \bigl(\langle
\beta\rangle+1\bigr)!.
\end{equation}
Selecting $b_0^*:= b_0 + \Delta$, $b_j^*:=b_j$, $j>0$, we
realize that
\begin{eqnarray}
\sum_{j=0}^{\beta^*} b_j^*
(x_i-x)^j & =& \sum_{j=0}^{\beta^*}
b_j (x_i-x)^j + \Delta\geq
Y_i\nonumber
\\
\eqntext{\forall i=1,\ldots,n\mbox{ with }\llvert
x-x_i\rrvert\leq h_k,}
\end{eqnarray}
so that by the definition of the $\hat{b}_j$, $j=0,\ldots,\beta
^*$, we have
%
\begin{eqnarray}
\label{eqlemcon2} \sum_{\llvert x_i-x\rrvert\leq h_k} \sum
_{j=0}^{\beta^*} \hat{b}_j
(x_i-x)^j & \leq&\sum_{\llvert x_i-x\rrvert\leq h_k}
\sum_{j=0}^{\beta^*} b_j^*
(x_i-x)^j
\nonumber\\[-8pt]\\[-8pt]\nonumber
&=& \sum_{\llvert x_i-x\rrvert\leq h_k}
\Biggl\{ \sum_{j=0}^{\beta
^*} b_j
(x_i-x)^j + \Delta\Biggr\}.
\end{eqnarray}
We define the polynomial
\[
Q(y):= \sum_{j=0}^{\beta^*} (
\hat{b}_j-b_j) (y-x)^j - \Delta.
\]
Then inequality (\ref{eqlemcon2}) implies that
%
\begin{equation}
\label{eqIntBed} \inf_{n,k,f} \int Q(y) \,d\lambda_n(y)
\leq0,
\end{equation}
where $\lambda_n$ denotes the uniform probability measure on the
discrete set $\{x_i\dvtx  \llvert x_i-x\rrvert\leq h_k\}$ inside the interval
$[x-h_k,x+h_k] \cap[0,1]$. We introduce the sets $Q^\pm$ of all $y\in
[x-h_k,x+h_k] \cap[0,1]$ such that $Q(y)$ is nonnegative or negative,
respectively. Our first task is to show that
%
\begin{equation}
\label{eqlemcon3} \inf_{n,k,f} \lambda_n\bigl(Q^-\bigr)
> 0\quad\mbox{or}\quad Q=0\qquad\mbox{identically}.
\end{equation}
In the latter case Theorem~\ref{teocon} is trivially true; hence we
focus on the case where $Q \neq0$. As $Q^-$ is the complement of $Q^+$
with respect to $[x-h_k,x+h_k] \cap[0,1]$, we have $\lambda
_n(Q^+)\geq1/2$ or $\lambda_n(Q^-)\geq1/2$. Clearly, we have (\ref
{eqlemcon3}) in the second case, so let us study the situation where
$\lambda_n(Q^+)\geq1/2$.

As $Q$ is a polynomial with degree $\leq\beta^*$ the set $Q^+$ equals
the union of at most $\beta^*+1$ disjoint sub-intervals of
$[x-h_k,x+h_k] \cap[0,1]$. The number of all design points in
$[x-h_k,x+h_k]$ is denoted by $m_k$. Hence, there exists at least one
interval $I_0^+ \subseteq Q^+$ such that $\lambda_n(I_0^+) \geq
1/(2\beta^*+2)$. At least $\lceil m_k / (2\beta^*+2) \rceil$ of the
$x_i$ lie in $I_0^+$ so that, due to the equidistant location of the
design points, the length of $I_0^+$ is larger or equal to
\[
\bigl\{\bigl\lceil m_k / \bigl(2\beta^*+2\bigr) \bigr\rceil- 1\bigr
\}/n \geq\bigl\{\bigl\lceil\lfloor n h_k \rfloor/ \bigl(2\beta^*+2
\bigr) \bigr\rceil- 1\bigr\}/n \geq c_1\bigl(\beta^*\bigr) \cdot
h_k, %
\]
for $n$ sufficiently large and some uniform constant $c_1(\beta^*)>0$
which does not depend on $n$ or $k$, but only on $\beta^*$. The
polynomial $Q$ takes\vspace*{1pt} only nonnegative values on the interval $I_0^+$.
By Lemma~\ref{lempoly} below there exists some interval $I_1^+
\subseteq I_0^+$ with the length $c_2(\beta^*) h_k$ such that
\[
\inf_{y\in I_1^+} \bigl\llvert Q(y)\bigr\rrvert\geq
c_3\bigl(\beta^*\bigr) \cdot\sup_{\llvert z-x\rrvert\leq h_k} \bigl
\llvert Q(z)\bigr\rrvert, %
\]
where the constants $c_2(\beta^*), c_3(\beta^*)>0$ only depend on
$\beta^*$. It follows from there that
\begin{eqnarray*}
\int_{Q^+} Q(y) \,d\lambda_n(y) &\geq&\int
_{I_1^+} Q(y) \,d\lambda_n(y)
\\
&\geq&
\lambda_n\bigl(I_1^+\bigr)\cdot\inf_{y\in I_1^+}
\bigl\llvert Q(y)\bigr\rrvert
\\
&\geq&\lambda_n\bigl(I_1^+
\bigr) c_3\bigl(\beta^*\bigr) \cdot\sup_{\llvert z-x\rrvert\leq h_k}
\bigl\llvert Q(z)\bigr\rrvert. %
\end{eqnarray*}
On the other hand we learn from (\ref{eqIntBed}) that
\[
\int_{Q^+} Q(y) \,d\lambda_n(y) \leq\int
_{Q^-} \bigl\llvert Q(y)\bigr\rrvert \,d\lambda_n(y)
\leq\lambda_n\bigl(Q^-\bigr)\cdot\sup_{\llvert z-x\rrvert\leq h_k}
\bigl
\llvert Q(z)\bigr\rrvert, %
\]
so that
\begin{eqnarray*}
\lambda_n\bigl(Q^-\bigr) &\geq&\lambda_n
\bigl(I_1^+\bigr) c_3\bigl(\beta^*\bigr) \geq
c_3\bigl(\beta^*\bigr) \cdot\bigl(c_2\bigl(\beta^*
\bigr) h_k n - 1\bigr) / m_k
\\
&\geq& c_3\bigl(
\beta^*\bigr) \cdot\bigl(c_2\bigl(\beta^*\bigr) - n^{-\mathfrak{h}_0}
\bigr) / \bigl(2 + n^{-\mathfrak{h}_0}\bigr), %
\end{eqnarray*}
unless $Q=0$ identically. Thus (\ref{eqlemcon3}) has been shown.

Using the arguments as above, we can now find some interval $I_0^-
\subseteq Q^-$ whose length is bounded from below by a constant (only
depending on $\beta^*$) times $h_k$. By Lemma~\ref{lempoly} there
exists an interval $I_1^- \subseteq I_0^-$, whose length is also
bounded from below by a constant (only depending on $\beta^*$) times
$h_k$ and on which $\llvert Q\rrvert$ is bounded from below by a
uniform multiple of
\[
\sup_{\llvert z-x\rrvert\leq h_k} \bigl\llvert Q(z)\bigr\rrvert\geq
\bigl\llvert
Q(x)\bigr\rrvert\geq\llvert\hat{b}_0 - b_0\rrvert-
\Delta. %
\]
This implies that
%
\begin{equation}
\label{eqlemcon4} \inf_{y \in I_1^-} \bigl(-Q(y)\bigr) \geq
c_4\bigl(\beta^*\bigr) \bigl(\llvert\hat{b}_0 -
b_0\rrvert- \Delta\bigr).
\end{equation}
On the other hand, for all $x_i \in I_1^-$ we have
%
\begin{eqnarray}
\label{eqlemcon5} Q(x_i) &=& \sum_{j=0}^{\beta^*}
\hat{b}_j (x_i-x)^j +
\Delta_i - f(x_i) - \Delta\geq Y_i -
f(x_i) - 2\Delta.
\end{eqnarray}
Combining the inequalities in (\ref{eqlemcon4}) and (\ref
{eqlemcon5}), we conclude that
\[
\bigl\llvert\tilde{f}_k(x) - f(x)\bigr\rrvert=\llvert
\hat{b}_0 - b_0\rrvert\leq- c^*\bigl(\beta^*\bigr)
\max\bigl\{\varepsilon_i\dvtx  x_i\in I_1^-
\bigr\} + c^*\bigl(\beta^*\bigr) \Delta%
\]
for some positive constant $c^*(\beta^*)$. Choosing $J(\beta
^*)$ sufficiently large (regardless of $k$, $n$ and $f$) there exists
some $l=1,\ldots,2J(\beta^*)$ such that $x+h_k {\mathcal I}_l
\subseteq I_1^-$, and hence
\[
\bigl\llvert\tilde{f}_k(x) - f(x)\bigr\rrvert\leq c^*\bigl(\beta^*
\bigr) \Delta- c^*\bigl(\beta^*\bigr)\cdot Z_j(h_k,x),
\]
which completes the proof.
\end{pf*}

%
\begin{lem} \label{lempoly}
Let $Q$ by any polynomial with the degree $\leq\beta^*$ and
$I\subseteq[x-h_k,x+h_k]$ be an interval with the length $\geq
c_5(\beta^*) h_k$ for some constant \mbox{$c_5(\beta^*)>0$}. Then there
exist some finite constants $c_6(\beta^*),c_7(\beta^*)>0$ which only
depend on $\beta^*$ and some interval $I^* \subseteq I$ with the
length $\geq c_6(\beta^*) h_k$ such that
\[
\inf_{y\in I^*} \bigl\llvert Q(y)\bigr\rrvert\geq
c_7\bigl(\beta^*\bigr) \cdot\sup_{\llvert z-x\rrvert\leq
h_k} \bigl
\llvert Q(z)\bigr\rrvert.
\]
\end{lem}

\begin{pf} 
If $Q$ is a constant function, the assertion is satisfied by putting
$c=1$. Otherwise, by the fundamental theorem of algebra, $Q$ can be
represented by
\[
Q(y) = \alpha_Q \prod_{j=1}^{\beta'}
(y-y_j), %
\]
where $1\leq\beta'\leq\beta^*$, the $y_j$ denote the
complex-valued roots of $Q$. By the pigeon hole principle there exists
some square $I_1^+ \times[-c_5(\beta^*) h_k/(2\beta^*+2),\break c_5(\beta
^*) h_k/(2\beta^*+2)]$ in the complex plane which does not contain any
$y_j$ where \mbox{$I_1^+ \subseteq I$} has the length $c_5(\beta^*)
h_k/(\beta^*+1)$. Now we shrink that square by the factor $1/2$ where
the center of the square does not change, leading to the square
$I_2^+ \times[-c_5(\beta^*) h_k/(4\beta^*+4),c_5(\beta^*)
h_k/(4\beta^*+4)]$. Thus, for any $y$ in this shrinked square, the
distance between $y$ and any $y_j$ is bounded from below by $c_5(\beta
^*) h_k/(4\beta^*+4)$ and by $\llvert y_j-x\rrvert- h_k$. If the
latter bound
dominates, we have $\llvert y_j - x\rrvert\geq\{1 + c_5(\beta^*)
/(4\beta^*+4)\}
\cdot h_k$. Then the distance between any $z \in[x-h_k,x+h_k]$ and
$y_j$ has the upper bound
\[
\llvert y_j-x\rrvert+ h_k \leq\llvert y_j-y
\rrvert+ 2h_k \leq\bigl\{1 + \bigl(8\beta^*+8\bigr)/c_5
\bigl(\beta^*\bigr)\bigr\} \cdot\llvert y_j-y\rrvert, %
\]
when applying the first bound. Otherwise, if the first bound dominates,
we have
\begin{eqnarray*}
\llvert z-y_j\rrvert&\leq&\llvert y_j-x\rrvert+
h_k \leq\bigl\{2 + c_5\bigl(\beta^*\bigr)/\bigl(4
\beta^*+4\bigr)\bigr\} \cdot h_k
\\
&\leq&\llvert y-y_j\rrvert
\cdot\bigl(4\beta^*+4\bigr) \bigl\{2/c_5\bigl(\beta^*\bigr) + 1/
\bigl(4\beta^*+4\bigr)\bigr\}. %
\end{eqnarray*}
In both cases $\llvert z-y_j\rrvert$ is bounded from above by a
uniform constant
$c_6(\beta^*)$ times $\llvert y-y_j\rrvert$. Then we learn from the
root-decomposition of the polynomial $Q$ that
\[
\inf_{y\in I^*} \bigl\llvert Q(y)\bigr\rrvert\geq
c_7\bigl(\beta^*\bigr) \cdot\sup_{\llvert z-x\rrvert\leq
h_k} \bigl
\llvert Q(z)\bigr\rrvert, %
\]
for some deterministic constant $c_7(\beta^*)>0$, which only depends
on $\beta^*$.
\end{pf}

\begin{pf*}{Proof of Proposition~\ref{propestimator}}
Part (a) follows directly from the definition of Lepski's method. For
part (b) we obtain from (\ref{eqcond1}), and repeated application of
the triangle and Jensen's inequality that
%
\begin{eqnarray}\label{eqlem00}
\quad&& \mathbb{E}_f \bigl[\llVert\hat{f} -
\tilde{f}_{\hat{k}^*}\rrVert^q \ind\bigl(\hat{k} <
\hat{k}^* \bigr) \bigr]^{1/q}\nonumber
\\
&&\qquad \leq\mathbb{E}_f \bigl[
\llVert\tilde{f}_{\hat{k}} - f\rrVert^q \ind\bigl(\hat{k}
< \hat{k}^* \bigr) \bigr]^{1/q} + \mathbb{E}_f \bigl[
\llVert\tilde{f}_{\hat{k}^*} - f\rrVert^q \ind\bigl(\hat{k}
< \hat{k}^* \bigr) \bigr]^{1/q}\nonumber
\\
&&\qquad \leq 2^{(q-2)/q}
\nonumber\\[-8pt]\\[-8pt]\nonumber
&&\quad\qquad{}\times \bigl( \mathbb{E}_f \bigl[
\bigl({R}_{\hat{k}}^q + B_{\hat{k}}^q\bigr)
\ind\bigl(\hat{k} < \hat{k}^* \bigr) \bigr]^{1/q} +
\mathbb{E}_f \bigl[\bigl({R}_{\hat{k}^*}^q +
B_{\hat{k}^*}^q\bigr) \ind\bigl(\hat{k} < \hat{k}^*
\bigr) \bigr]^{1/q} \bigr)
\\
&&\qquad  \leq 2^{(2q-1)/{q}} \bigl(\mathbb{E}_f \bigl[
B_{\hat{k}^*}^q \bigr]^{1/q} + \mathbb{E}_f
\bigl[{R}_{\hat{k}}^q \ind\bigl(\hat{k} < \hat{k}^*
\bigr) \bigr]^{1/q} \bigr)\nonumber
\\
&&\qquad  \leq 2^{(2q-1)/q} \Biggl( \mathbb{E}_f \bigl[
\hat{\z}_{\hat{k}^*}^q\bigr]^{1/q} + \sum
_{k=0}^{K-1} \mathbb{E}_f
\bigl[{R}_k^q \ind\bigl(\hat{k} = k, k <
\hat{k}^* \bigr) \bigr]^{1/q} \Biggr),\nonumber
\end{eqnarray}
where we also used that $R_k$ decreases in $k$ and $B_k$ increases in $k$.
Note that
\begin{eqnarray*}
\ind\bigl(\hat{k} = k, k < \hat{k}^* \bigr) & \leq&\ind\bigl
(\exists l
\leq k\dvtx  \llVert\tilde{f}_{k+1} - \tilde{f}_l\rrVert>
\hat{\z}_l^T + \hat{\z}_{k+1}^T
\bigr)\cdot\ind\bigl(k < \hat{k}^* \bigr)
\\
& \leq&\ind\bigl(\exists l \leq k\dvtx  \llVert\tilde{f}_{k+1} - f\rrVert
+ \llVert f - \tilde{f}_l\rrVert> \hat{\z}_l^T
+ \hat{\z}_{k+1}^T \bigr) \cdot\ind\bigl(k <
\hat{k}^* \bigr)
\\
\nonumber
& \leq&\ind\bigl(\exists l \leq k\dvtx  {R}_l +
B_{k+1} > \hat{\z}_l^T \bigr) \cdot\ind
\bigl(B_{k+1} \leq\hat{\z}_{k+1}^T/2 \bigr)
\\
& \leq&\ind\bigl(\exists l \leq k\dvtx  {R}_l > \hat{\z}_l^T/2 \bigr).
\end{eqnarray*}
Inserting this inequality into (\ref{eqlem00}) completes the proof.
\end{pf*}

In the sequel, the following three lemmas will be useful. The proofs
are given in the supplementary material~\cite{jirmeireisssuppl}.
%
\begin{lem}\label{lemquantcoomp}
If $y,t \to\infty$ and
\begin{eqnarray*}
y &=& \cf(\log t)^{\bd} t^{\aexp} \bigl(1 + \mbox{\scriptsize$
\mathcal{O}$}(1) \bigr), \qquad\cf, \aexp> 0, \bd\in\mathbb{R},
\end{eqnarray*}
then
\begin{eqnarray*}
t &=& \bigl(\cf^{-1} \bigl(\log y^{1/\aexp}\bigr)^{-\bd} y
\bigr)^{1/\aexp} + \OO(1).
\end{eqnarray*}
In particular, if we have $v = \mathcal{U}_x(y)$ with $v \to0$, then
\begin{eqnarray*}
F_x (v ) &=& 1 - \cf_x^{-\aexp_x} \bigl(\log
\llvert v\rrvert^{-1/\aexp
_x} \bigr)^{-\bd_x \aexp_x}\llvert v\rrvert
^{\aexp_x} \bigl(1 + \mbox{\scriptsize$\mathcal{O}$}(1) \bigr).
\end{eqnarray*}
\end{lem}

%
\begin{lem}\label{lemcomputeprodprob}
For $1 \leq j_0,j_1 \leq n$, let $\mathcal{J} = \{j_0,\ldots,j_1 \}$
such that $\llvert j_0 - j_1\rrvert/n = \OO(n^{-\rho_0}
)$ for some $0 < \rho_0 < 1$. If $u \to0$, $u \leq-n^{-\rho_1}$ for
some $\rho_1 > 0$, then
\begin{eqnarray*}
\prod_{j \in\mathcal{J}} P \bigl(\varepsilon_j \leq
A_{x_{j_0}}\bigl(-u^{-1}\bigr) \bigr) & \leq& e^{\# \mathcal{J} c_3^- u},
\end{eqnarray*}
where $c_3^- < 1$ may be chosen arbitrarily close to one.
\end{lem}

%
\begin{lem} \label{lempowern}
Let $(q_n)_n$ be a real-valued sequence which satisfies $q_n \in
[1,\log n]$ for all integer $n$, and denote with $F(\cdot)$ the c.d.f. of
$\varepsilon$. Then we have
\[
\mathbb{E}\bigl\llvert\max\{\varepsilon_1,\ldots,
\varepsilon_n\}\bigr\rrvert^{q_n} \leq\bigl(1 + \mbox{
\scriptsize$\mathcal{O}$}(1) \bigr)\int_0^{n^{1/2}}
\bigl((-\mathcal{U})^{q_n} (n/y ) \bigr)^{(1)} \exp(-y)\,dy.
\]
If $\mathcal{U}(\cdot)$ is not differentiable, replace $\mathcal{U}(\cdot)$
with $c_2^+ A(\cdot)$ in the above inequality, where~$c_2^+ > 1$ can be
chosen arbitrarily close to one. 
If $q_n = q$ is finite and independent of $n$, we obtain that
\[
\mathbb{E}\bigl\llvert\max\{\varepsilon_1,\ldots,
\varepsilon_n\}\bigr\rrvert^{q} = \OO\bigl((\log
n)^{q\bd_F} n^{-q/\aexp_F} \bigr). %
\]
For arbitrary $q_n \in[1,\log n]$ we have
\begin{eqnarray*}
\OO\bigl(n^{-c_2^+\aexp_F/q_n} \bigr) &\leq&\int_0^{\infty}
 F\bigl(-x^{1/{q_n}}\bigr)^n \,dx \leq\OO\bigl(n^{-c_2^-\aexp_F/q_n}
\bigr),
\end{eqnarray*}
where $0 < c_2^- < 1 < c_2^+$ can be chosen arbitrarily close to one.
\end{lem}

\begin{pf*}{Proof of Theorem~\ref{Tpointwise}}
In the course of the proof we will frequently apply Proposition~\ref
{propestimateunivquantile}. We may do so since condition $\mathfrak
{h}_0 < \beta_{\zz} \aexp_{\zz} / (\beta_{\zz} \aexp_{\zz} +
1)$ implies (\ref{eqparamrelation1}). The general strategy is the
following. By the triangle inequality and Jensen's inequality, we have
%
\begin{equation}
\quad\mathbb{E}_f \bigl[ \bigl(\hat{f}(\zz) - f(\zz)
\bigr)^2 \bigr] \leq2 \mathbb{E}_f \bigl[ \bigl(
\hat{f}(\zz) - \tilde{f}_{\hat{k}^*} \bigr)^2 \bigr] + 2
\mathbb{E}_f \bigl[ \bigl(\tilde{f}_{\hat{k}^*}- f(\zz)
\bigr)^2 \bigr],
\end{equation}
and we will treat both quantities separately. In order to deal with the
first, Proposition~\ref{propestimator} implies that it suffices to consider
\begin{eqnarray*}
&& \sup_{f \in H_{\NN(\zz)}(\beta,L)} \mathbb{E}_f \bigl[\bigl(\hat{\z
}_{\hat{k}^*}^T\bigr)^2 \bigr]^{1/2} +
\sup_{f \in H_{\NN(\zz
)}(\beta,L)} \sum_{k=0}^{K-1}
\mathbb{E}_f \bigl[{R}_k^{2} \ind\bigl(
\exists l \leq k\dvtx  {R}_l > \hat{\z}_l^T/2
\bigr) \bigr]^{1/2}
\\
&&\qquad =:I + \sum_{k=0}^{K-1} \mathit{II}_k.
\end{eqnarray*}
To treat $I$, we require the following simple lemma; the proof is given
in the supplementary material~\cite{jirmeireisssuppl}.
%
\begin{lem}\label{lemdealwithz}
Let $q \geq1$. Under the assumptions of Theorem~\ref{Tpointwise}, we
have uniformly over $f \in H_{\NN(\zz)}(\beta,L)$
\begin{eqnarray*}
\mathbb{E}_f \bigl[\bigl(\hat{\z}_{\hat{k}^*}^T
\bigr)^q \bigr] & \leq&\bigl(c_1^+\bigr)^q
\mathbb{E}_f \bigl[\bigl(\z_{\hat{k}^*}^T
\bigr)^q \bigr] + \OO\bigl(n^{-q/\aexp_{\zz}} \bigr),
\end{eqnarray*}
where $c_1^+>1$.
\end{lem}
Applying the above result with $q = 2$ we obtain
%
\begin{equation}
I^2 \leq\bigl(c_1^+\bigr)^2
\mathbb{E}_f \bigl[\bigl(\z_{\hat{k}^*}^T
\bigr)^2 \bigr] + \OO\bigl(n^{-2/\aexp_{\zz}} \bigr),
\end{equation}
and it remains to deal with $\mathbb{E}_f [(\z_{\hat{k}^*}^T)^2 ]$. We define
\[
k^\pm:= \inf\bigl\{k=0,\ldots,K-1\dvtx  B_{k+1} >
c_2^\pm\z_{k+1}^T / 2\bigr\} \wedge
K. %
\]
On the event ${\mathcal{A}}_n=\{c_2^- \z_k \leq\hat{\z}_k \leq
c_2^+ \z_k\mbox{ for all } k=0,\ldots,K-1\}$ we have $k^- \leq
\hat{k}^* \leq k^+$. From Proposition~\ref
{propestimateunivquantile} we infer $P(\mathcal{A}_n^c)=\OO
(n^{-v} \log n )$ due to\break \mbox{$K = \OO(\log n )$}. Since $\z
_k$ decreases monotonically in $k$, we deduce $\mathbb{E} [(\z
_{\hat{k}^*}^T)^2 ] \leq\z_{k^-}^2 + \OO
(n^{-2/\aexp_{\zz}} )$. Note that the deterministic sequences
$(h_k^{\pm})$ satisfy (see Lemma~\ref{lemquantcoomp} for details)
%
\begin{eqnarray}\label{eqthmpointwise15}
\z_{k^-} &\thicksim&- \bigl(h_k^{\pm}
\bigr)^{\beta_{\zz}} \quad\mbox{and}
\nonumber\\[-8pt]\\[-8pt]\nonumber
h_{k^{\pm}} &\thicksim& (n /\log n)^{-1/(\aexp
_{\zz} \beta_{\zz}+1)} (\log n)^{(\aexp_{\zz} \bd_{\zz})/(
\aexp_{\zz} \beta_{\zz} + 1)},
\end{eqnarray}
under our assumption $\mathfrak{h}_0 < \beta_0 \aexp_0 / (\beta_0
\aexp_0 + 1)$. We obtain
%
\begin{eqnarray}\label{eqthmpointwise2}
&& \sup_{f\in H_{\NN(\zz)}(\beta,L)} \mathbb{E}_f
\bigl[
\bigl(\hat{\z}{}^T_{\hat{k}^*}\bigr)^2\bigr]
\nonumber\\[-8pt]\\[-8pt]\nonumber
&&\qquad = \OO
\bigl((n /\log n)^{(-2\beta_{\zz})/(\aexp_{\zz} \beta_{\zz}+1)} (\log n)^{(2\aexp
_{\zz} \bd_{\zz}\beta_{\zz})/(\aexp_{\zz} \beta_{\zz}+1)} \bigr),
\end{eqnarray}
and it remains to deal with the second part. Let $\mathcal{B}_k =
\{\exists l \leq k\dvtx  {R}_l > \hat{\z}_l^T/2 \}$.
Then for $\delta_0 = n^{-\mathfrak{h}_0/\aexp_{\zz}}$, we have
\begin{eqnarray*}
\mathit{II}_k^2 &=& \int_0^{\infty}P_f
\bigl(R_k^2 \ind(\mathcal{B}_k)\geq x \bigr)
\,dx \leq\int_0^{\delta_0} P_f (
\mathcal{B}_k )\,dx + \int_{\delta_0}^{\infty}P
\bigl(R_k^2 \geq x \bigr)\,dx
\\
&\leq&\delta_0 \sum_{l = 0}^k
P_f \bigl({R}_l > \hat{\z}_l^T/2
\bigr) + \int_{\delta_0}^{\infty}P \bigl(R_k^2
\geq x \bigr)\,dx =: \mathit{III}_k + \mathit{IV}_k.
\end{eqnarray*}
We first deal with $\mathit{III}_k$. Recall that for $j=1,\ldots,2J(\beta^*)$
we have
\[
Z_j(h_k,\zz) = \max\{\varepsilon_i\dvtx
x_i \in\zz+ h_k {\mathcal I}_j \}
\]
and
\[
{\mathcal I}_j:= \bigl[-1+(j-1)/J\bigl(\beta^*
\bigr),-1+j/J\bigl(\beta^*\bigr)\bigr].
\]
Put\vspace*{1pt} $\mathcal{J}_j(l) = \{0,1/n,2/n,\ldots,1 \}\cap
\{\zz+ h_l {\mathcal I}_j \}$ for $j=1,\ldots,2J(\beta^*)$, and
note that $\#\mathcal{J}_j(l) \geq n h_l/2 J(\beta^*)$. An
application of Proposition~\ref{propestimateunivquantile} and Theorem~\ref{teocon} then yields that
%
\begin{eqnarray}
P_f \bigl({R}_l > \hat{\z}_l^T/2 \bigr) &\leq&\sum_{j =
1}^{J(\beta^*)}P
\bigl( \bigl\llvert Z_j(h_l,\zz)\bigr\rrvert>
c_1^- \z_l \bigr) + \OO\bigl(n^{-2/\aexp_{\zz}} \bigr)
\nonumber\\[-8pt]\\[-8pt]\nonumber
&\leq&\sum_{j = 1}^{2J(\beta
^*)} \prod
_{i \in\mathcal{J}_j(l)} P \bigl(\varepsilon_i <
-c_1^-\z_l \bigr) + \OO\bigl(n^{-2/\aexp_{\zz}}
\bigr),
\end{eqnarray}
where $c_1^- <1$ may be chosen arbitrarily close to one. Arguing as in
Lemma~\ref{lemcomputeprodprob} 
we obtain
%
\begin{equation}
\label{eqthmpointwise3} \prod_{i \in\mathcal{J}_j} P \bigl(
\varepsilon_i > -c_1^-\z_l \bigr) = \OO
\bigl(n^{-2c_2^-/\aexp_{\zz}} \bigr),
\end{equation}
where $c_2^- <1$ may be chosen arbitrarily close to one. Hence we
obtain that
%
\begin{equation}
\label{eqthmpointwise4} \qquad \mathit{III}_k \leq\delta_0 \sum
_{l = 0}^k \sum_{j = 1}^{2J(\beta^*)}
\OO\bigl(n^{-2c_2^-/\aexp_{\zz}} \bigr) = \OO\bigl(K \delta_0
n^{-2c_2^-/\aexp_{\zz}} \bigr) = \OO\bigl(n^{-2/\aexp_{\zz}} \bigr),
\end{equation}
since\vspace*{1pt} $K = \OO(\log n )$. For dealing with $\mathit{IV}_k$, set
$\eta_0 = \exp( n^{\mathfrak{h}_0/4} )$. Let $u_{\zz} =
A_{\zz}^{-1}(\delta_0)$. Then Lemma~\ref{lemquantcoomp} implies that
$u_x < -n^{-\mathfrak{h}_0/2}$ for large enough $n$. Then as in (\ref
{eqthmpointwise3}), it follows from Lemma~\ref{lemcomputeprodprob}
that for sufficiently large $n$, 
\begin{eqnarray*}
\int_{\delta_0}^{\eta_k}P \bigl(R_k^2
\geq x \bigr)\,dx &=& \int_{\delta_0}^{\eta_k}P
\bigl(R_k \leq-x^{1/2} \bigr)\,dx
\\
&\leq&\eta_k \sum
_{j = 1}^{2J(\beta^*)} \prod
_{i \in\mathcal{J}_j(k)} P \bigl(\varepsilon_i < A_{\zz}
\bigl(-n^{-\mathfrak{h}_0/2}\bigr) \bigr)
\\
&\leq& \eta_k \sum_{j = 1}^{2J(\beta^*)}
\exp\bigl( -\#\mathcal{J}_{j(k)} n^{-\mathfrak{h}_0/2} \bigr)
\\
&\leq&
\eta_k \sum_{j =
1}^{2J(\beta^*)} \exp
\bigl( -n^{\mathfrak{h}_0/2} \bigr) = \OO\bigl(\exp\bigl(
-n^{\mathfrak{h}_0/4} \bigr)
\bigr).
\end{eqnarray*}
%
Let $p> 2$. Since $\mathbb{E} [\llvert R_k\rrvert^p ] = \OO(1
)$ by Lemma~\ref{lempowern}, it follows from the Markov inequality that
\begin{eqnarray*}
\int_{\eta_k}^{\infty}P \bigl(R_k^2
\geq x \bigr)\,dx & \leq&\int_{\eta_k}^{\infty}x^{-p/2}
\mathbb{E} \bigl[\llvert R_k\rrvert^p \bigr] \,dx =
\frac{2}{p - 2} \eta_k^{-p/2 + 1} \OO(1 ).
\end{eqnarray*}
Combining the above and (\ref{eqthmpointwise4}), it follows that $\mathit{IV}_k
= \OO(\exp( -n^{\mathfrak{h}_0/8} ) )$, which
in turn yields
%
\begin{equation}
\mathit{II}_k^2 = \mathit{III}_k + \mathit{IV}_k = \OO
\bigl(n^{-2/\aexp_{\zz}} \bigr) + \OO\bigl(\exp\bigl( -n^{\mathfrak
{h}_0/8} \bigr)
\bigr) = \OO\bigl(n^{-2/\aexp_{\zz}} \bigr).
\end{equation}
We thus conclude
\begin{eqnarray*}
\sum_{k=0}^{K-1} \mathit{II}_k &=& \OO
\bigl(K n^{-1/\aexp_{\zz}} \bigr) = \OO\bigl(\log n n^{-1/\aexp_{\zz}}
\bigr)
\\
&=&
\mbox{\scriptsize$\mathcal{O}$} \bigl((n \log n)^{- \beta_{\zz
}/(\aexp_{\zz}
\beta_{\zz}+1)} (\log
n)^{(\aexp_{\zz} \bd_{\zz})/({\aexp
_{\zz} \beta_{\zz}+1})} \bigr).
\end{eqnarray*}
Piecing everything together and taking squares, we arrive at
%
\begin{eqnarray}\label{eqthmpointwise5}
&& \sup_{f \in H_{\NN(\zz)}(\beta,L)}\mathbb{E}_f
\bigl[
\bigl(\hat{f}(\zz) - \tilde{f}_{\hat{k}^*} \bigr)^2 \bigr]
\nonumber\\[-8pt]\\[-8pt]\nonumber
&&\qquad =
\OO\bigl((n \log n)^{- \beta_{\zz}/({\aexp_{\zz} \beta_{\zz
}+1})} (\log n)^{(\aexp_{\zz} \bd_{\zz})/(\aexp_{\zz} \beta
_{\zz}+1)} \bigr).
\end{eqnarray}
To complete the proof, it remains to deal with $\mathbb{E}_f
[ (\tilde{f}_{\hat{k}^*} - f(\zz) )^2 ]$. Let $p'
> 2$. Then by (\ref{EqRk2}) and the triangle, Jensen and H\"{o}lder
inequalities, we have
%
\begin{eqnarray}
\mathbb{E}_f \bigl[ \bigl(\hat{f}(\zz) -
\tilde{f}_{\hat{k}^*} \bigr)^2 \bigr] &\leq& 2 \mathbb{E}_f
\bigl[R_{\hat{k}^*}^2 + B_{\hat{k}^*}^2 \bigr]
\nonumber\\[-8pt]\label{eqthmpointwise6}  \\[-8pt]\nonumber
 & \leq&2 \mathbb{E}_f \bigl[R_{\hat{k}^*}^2\ind(
\mathcal{A}_n ) + B_{\hat{k}^*}^2\ind(
\mathcal{A}_n ) \bigr]
\\
&&{}+ P \bigl(\mathcal{A}_n^{c} \bigr)^{(p'-2)/p'}
\bigl(\mathbb{E} \bigl[R_{k^-}^{p'} \bigr]^{1/p'} +
\OO(L_{\zz} ) \bigr).
\end{eqnarray}
Hence Proposition~\ref{propestimateunivquantile}, Lemma~\ref{lempowern} and (\ref{eqthmpointwise15}) imply that
the above is of order
%
\begin{eqnarray}
\label{eqthmpointwise7} \mathbb{E}_f \bigl[R_{\hat{k}^*}^2
+ B_{\hat{k}^*}^2 \bigr] &\leq&2 \mathbb{E}_f
\bigl[R_{{k}^-}^2 \bigr] + B_{k^+}^2 +
\OO\bigl(n^{-2/\aexp_{\zz}} \bigr)\nonumber
\\
&=&
\OO\bigl(n_{k^-}^{-2/\aexp
_{\zz}} (\log
n_{k^-})^{2 \bd_{\zz}} + \z_{k^+}^2 +
n^{-2/\aexp
_{\zz}} \bigr)
\\
&=& \OO\bigl((n \log n)^{(- 2\beta_{\zz
})/(\aexp_{\zz} \beta_{\zz}+1)} (\log n)^{(2\aexp_{\zz} \bd
_{\zz})/(\aexp_{\zz} \beta_{\zz}+1)} \bigr).\nonumber
\end{eqnarray}
The above bound is uniform over $H_{\NN(\zz)}(\beta,L)$, and the
proof is complete.
\end{pf*}

\begin{pf*}{Proof of Theorem~\ref{TL2upperbound}}
For the proof we require the following lemma, which provides a
sub-polynomial upper bound on the probability that $R_k$ exceeds the
threshold $\hat{\z}_l^T/2$.
%
\begin{lem} \label{LL2estimate}
Suppose $h_k \leq\exp(-c_H \log^\gamma n)$ for fixed constants
$\gamma\in(0,1)$, $c_H>0$. Grant Assumption~\ref{assmain}, and let
$\mathfrak{m}, \aexp_0, \beta_0, \mathfrak{h}_0$ satisfy (\ref
{eqparamrelation1}) in view of (\ref{eqlowbnd}). Then
\[
\sup_{k = 0,\ldots,K-1} P \bigl(R_k > \hat{\z}_l^T/2 \bigr) = \OO\bigl
(\exp\bigl(-c_H
\log^{1+\gamma} n / 2q\bigr) \bigr),
\]
as $n\to\infty$.
\end{lem}
Let $p> q$. Then Lemma~\ref{lempowern} gives
%
\begin{equation}
\label{eqRP} \mathbb{E} \bigl[R_k^p
\bigr]^{q/p} = \OO\bigl((nh_k)^{-q /\aexp_F} (\log n
h_k)^{q \bd_F} \bigr).
\end{equation}
%
An application of H\"{o}lder's inequality, (\ref{eqRP}) and Lemma~\ref{LL2estimate} yields that
%
\begin{eqnarray} \label{eqest1}
&& \sup_{f\in H_{[0,1]}(\beta,L)} \sum_{k=0}^{K-1}
\mathbb{E}_f \bigl[{R}_k^{q} \ind\bigl(
\exists l \leq k\dvtx  {R}_l > \hat{\z}_l^T/2
\bigr) \bigr]^{1/q}\nonumber
\\
&&\qquad \leq\sum_{k=0}^{K}
(k+1)^{(p-q)/p}1_{[0,\exp(-c_H \log^\gamma n)]}(h_k)\nonumber
\\
&&\hspace*{45pt}{}\times 
\OO \bigl(n^{-\mathfrak{h}_0 /\aexp_F}\cdot\exp\bigl(-c_H/q (p-q)
\log^{1+\gamma} n / [2 q p] \bigr) \bigr)
\\
&&\quad\qquad{} + \sum_{k=0}^{K} 1_{(\exp(-c_H \log^\gamma n),\infty)}(h_k)
\cdot(nh_k)^{-1/\aexp_F}\nonumber
\\
&&\qquad  = \OO\bigl(n^{- c_2^-/\aexp_F} \exp
\bigl([c_H/\aexp_F] \log^\gamma n\bigr)\cdot
\log n \bigr).\nonumber
\end{eqnarray}
Choosing $c_H, \gamma> 0$ sufficiently small, the above is of order
$\OO(n^{-c_3^-/\aexp_F} )$, where $c_3^- < 1$ can be chosen
arbitrarily close to one. According to Proposition~\ref{propestimator}
it remains to bound the expectation $\mathbb{E}_f [(\hat{\z}{}^T_{\hat
{k}^*})^q]$ uniformly over $f\in H_{[0,1]}(\beta,L)$.
Applying Lemma~\ref{lemdealwithz}, we obtain that uniformly over $f\in
H_{[0,1]}(\beta,L)$
\begin{eqnarray*}
\mathbb{E}\bigl[\bigl(\hat{\z}_{\hat{k}^*}^T
\bigr)^q\bigr] & \leq&\bigl(T c_1^+\bigr)^q
\mathbb{E}_f \bigl[(\z_{\hat{k}^*})^q\bigr] + \OO
\bigl(n^{-q/\aexp_{F}} \bigr).
\end{eqnarray*}
To deal with $\mathbb{E}_f [(\z_{\hat{k}^*})^q]$, we introduce
\[
k^\pm:= \inf\bigl\{k=0,\ldots,K-1\dvtx  B_{k+1} >
c_2^\pm\z_{k+1}^T / 2\bigr\} \wedge
K. %
\]
On the event $\mathcal{A}_n=\{c_2^- \z_k \leq\hat{\z}_k \leq
c_2^+ \z_k\mbox{ for all } k=0,\ldots,K-1\}$ we have $k^- \leq
\hat{k}^* \leq k^+$. From Proposition~\ref
{propestimateunivquantile} we infer $P(\mathcal{A}_n^c)=\OO
(n^{-q/\aexp_F} )$. Since $\z_k$ decreases monotonically in $k$,
we find $\mathbb{E}[\z_{\hat{k}^*}^q] \leq\z_{k^-}^q +
\OO(n^{-q/\aexp_F} )$. Note that the deterministic sequences
$(h_k^{\pm})$ satisfy
%
\begin{eqnarray}
h_{k^{\pm}} &\thicksim& (n)^{-1/(\aexp_{F} \beta_+1)} (\log
n)^{(\aexp_{F} \bd_{F})/(\aexp_{F} \beta_+1)} \quad\mbox{and}
\nonumber\\[-8pt]\\[-8pt]\nonumber
\z_{k^-} &\thicksim& n^{-\beta/(\aexp_{F} \beta
_+1)} (\log n)^{(\beta\aexp_{F} \bd_{F})/(\aexp_{F} \beta_+1)},
\end{eqnarray}
provided that $\mathfrak{h}_0 < \beta_0 \aexp_0 / (\beta_0 \aexp_0
+ 1)$. For computational details, refer to Lem\-ma~\ref{lemquantcoomp}.
Moreover, condition $\mathfrak{h}_0 < \beta_0 \aexp_0 / (\beta_0
\aexp_0 + 1)$ also implies (\ref{eqparamrelation1}). We obtain
%
\begin{equation}
\label{eqest2} \quad\sup_{f\in H_{[0,1]}(\beta,L)} \mathbb{E}_f \bigl[
\bigl(\hat{\z}{}^T_{\hat{k}^*}\bigr)^q \bigr] = \OO
\bigl(n^{(-q\beta)/(\aexp_{F} \beta_+1)} (\log n)^{(q\beta\aexp_{F} \bd
_{F})/(\aexp_{F} \beta_+1)} \bigr).
\end{equation}
Combining this result with (\ref{eqest1}), Proposition~\ref
{propestimator} yields that
\[
\sup_{f\in H_{[0,1]}(\beta,L)} \mathbb{E}_f \bigl[ \llVert
\hat{f} - \tilde{f}_{\hat{k}^*}\rrVert_q^q \bigr]
= \OO\bigl(n^{(-q\beta)/(\aexp_{F} \beta_+1)} (\log n)^{(q\beta\aexp_{F}
\bd_{F})/(\aexp_{F} \beta_+1)} \bigr).
\]
Arguing similarly as in (\ref{eqthmpointwise7}), by (\ref{EqBk}) and
(\ref{EqRk}) we deduce that
%
\begin{eqnarray}\label{eqoracleestimator}
\mathbb{E}_f \bigl[ \llVert\tilde{f}_{\hat{k}^*} -
f\rrVert_q^q \bigr] & \leq&2^{q}
\mathbb{E}_f \bigl[B_{\hat{k}^*}^q\bigr] +
2^{q} \mathbb{E}_f \bigl[R_{\hat{k}^*}^q
\bigr]\nonumber
\\
&\leq&2^{q} B_{k^+}^q + \OO
\bigl(n^{-q/\aexp_F} \bigr) + 2^{q} \mathbb{E}_f
\bigl[R_{k^-}^q\bigr]
\\
&=&\OO\bigl(n^{(-q\beta)/(\aexp_{F} \beta
_+1)} (\log n)^{(q\beta\aexp_{F} \bd_{F})/(\aexp_{F} \beta_+1)}
\bigr),\nonumber
\end{eqnarray}
uniformly with respect to $f\in H_{[0,1]}(\beta,L)$, by conditioning
on the event ${\mathcal{A}}_n$ and using (\ref{eqRP}). The proof is complete.
\end{pf*}

The proofs of Theorems~\ref{Tlowerbound},~\ref{T421} and~\ref{teo-aexp2} are given in the supplementary material~\cite{jirmeireisssuppl}.

\section*{Acknowledgments}
We would like to thank the anonymous referees for many
constructive remarks that have lead to a significant improvement both
in the results and the presentation. We also thank Holger Drees and
Keith Knight for extensive discussions and insightful comments on
quantile estimation in inhomogeneous data.


\begin{supplement}[id=suppA]\label{suppA}
\stitle{Additional simulations, proof of lower bound, technical lemmas and sharpness estimation}
\slink[doi]{10.1214/14-AOS1248SUPP} 
\sdatatype{.pdf}
\sfilename{aos1248\_supp.pdf}
\sdescription{In the supplementary material we provide additional simulations and the proofs
of the lower bound results as well as technical lemmas and sharpness
estimation.}
\end{supplement}

%

\printaddresses
\end{document}